\documentclass[a4paper]{amsart}
\usepackage[utf8]{inputenc}
\usepackage{amsfonts}
\usepackage{amsmath,latexsym,amssymb,amsfonts, amsthm}
\usepackage{esint}
\usepackage{graphicx}
\usepackage{amscd}
\usepackage{color}
\usepackage{bm}             
\usepackage{enumerate}
\usepackage{aliascnt}
\usepackage[dvips]{epsfig}
\usepackage{psfrag}
\usepackage{enumitem}
\usepackage{epsfig}
\usepackage{comment}

\usepackage[
  hmarginratio={1:1},     
  vmarginratio={1:1},     
  textwidth=16cm,        
  textheight=21cm,
  heightrounded,          
]{geometry}

\usepackage{graphicx,color}
\usepackage[colorlinks]{hyperref}
\usepackage{cleveref}
\hypersetup{linkcolor=blue,citecolor=blue,filecolor=black,urlcolor=blue}
\usepackage{mathtools}

\DeclarePairedDelimiter{\norm}{\lVert}{\rVert}

\renewcommand{\b}{\beta}

\definecolor{my_red}{RGB}{255, 0, 0}
\newcommand{\red}[1]{{\color{my_red}#1}}

\newtheorem{theorem}{Theorem}[section]
\crefname{theorem}{theorem}{theorems}
\Crefname{theorem}{Theorem}{Theorems}

\newaliascnt{proposition}{theorem}
\newtheorem{proposition}[proposition]{Proposition}
\aliascntresetthe{proposition}
\crefname{proposition}{proposition}{propositions}
\Crefname{proposition}{Proposition}{Propositions}

\newaliascnt{corollary}{theorem}

\aliascntresetthe{corollary}
\crefname{corollary}{corollary}{corollaries}
\Crefname{corollary}{Corollary}{Corollaries}

\newaliascnt{lemma}{theorem}
\newtheorem{lemma}[lemma]{Lemma}
\aliascntresetthe{lemma}
\crefname{lemma}{lemma}{lemmas}
\Crefname{lemma}{Lemma}{Lemmas}

\theoremstyle{definition}
\newaliascnt{remark}{theorem}
\newtheorem{remark}[remark]{Remark}
\aliascntresetthe{remark}
\theoremstyle{definition}
\crefname{remark}{remark}{remarks}
\Crefname{remark}{Remark}{Remarks}

\newaliascnt{example}{theorem}

\aliascntresetthe{example}
\theoremstyle{definition}
\crefname{example}{example}{examples}
\Crefname{example}{Example}{Examples}

\newaliascnt{definition}{theorem}

\aliascntresetthe{definition}
\theoremstyle{definition}
\crefname{definition}{definition}{definitions}
\Crefname{definition}{Definition}{Definitions}

\numberwithin{theorem}{section}

\numberwithin{equation}{section}

\newcommand{\rad}{\mathrm{rad}}

\newcommand{\beq}{\begin{equation}}
\newcommand{\eeq}{\end{equation}}

\newcommand{\R}{\mathbb{R}}
\newcommand{\N}{\mathbb{N}}

\newcommand{\mf}[1]{\mathbf{#1}}
\newcommand{\bs}[1]{\boldsymbol{#1}}

\newcommand{\pa}{\partial}

\newcommand{\cA}{{\mathcal A}}   
\newcommand{\cB}{{\mathcal B}}

\newcommand{\cM}{{\mathcal M}}   
\newcommand{\cN}{{\mathcal N}}

\newcommand{\cS}{{\mathcal S}}

\newcommand{\intRn}{\int_{\R^d}}
\usepackage{mathtools}

\DeclareMathOperator*{\argmin}{arg\,min}

\newcommand{\dist}{{\mathrm{dist}}}
\newcommand{\supp}{{\mathrm{supp}}}

\newcommand{\eps}{\varepsilon}

\usepackage{accents}

\renewcommand{\epsilon}{\varepsilon}

\author[H. Fan, L. Giaretto, N. Soave]{Haining Fan, Lorenzo Giaretto and Nicola Soave}\thanks{}
\address{Haining Fan \newline \indent School of mathematics, China University of Mining and Technology,\newline \indent Xuzhou, Jiangsu 221116, China}
\email{fanhaining888@163.com}
\address{Lorenzo Giaretto and Nicola Soave\newline \indent
 Dipartimento di Matematica ``Giuseppe Peano'', Universit\`a di Torino, \newline \indent
Via Carlo Alberto 10,
10123 Torino, Italy}
\email{lorenzo.giaretto@unito.it, nicola.soave@unito.it}

\title[On least energy solutions for a critical NLS system with $K$-wise interaction]{On least energy solutions for a critical nonlinear Schr\"odinger system with $K$-wise interaction: the competitive case
}
\thanks{L.G. and N.S. are affiliated to the INDAM - GNAMPA group. L.G. was partially supported by the INDAM--GNAMPA 2026 research project ``Ottimizzazione Geometrica e Frontiere Libere''.
H.F. is supported by the Fundamental Research Funds for the Central Universities (No. 2024KYJD2007)
and NSFC (No.12671144).\\
Declarations of interest: none. \\
Data availability: Data sharing not applicable to this article as no datasets were generated or analyzed during the current study.}

\begin{document}

\begin{abstract}
We study least energy fully non-trivial solutions of a nonlinear
Schr\"odinger system with Sobolev critical nonlinearities and repulsive interaction involving all the components at the same time.
Under the scalar Brezis-Nirenberg assumptions on the linear coefficients,
we prove existence of least energy fully non-trivial
solutions on smooth bounded domains for every negative coupling
in dimensions \(N=5,6\).
For \(N=3,4\), due to the lack of $L^2$-integrability of the Aubin-Talenti bubbles, the same conclusion holds under additional sufficient assumptions, and in particular we provide an explicit
condition on the domain and coefficients.
We also investigate the strong competition limit, and characterize fully non-trivial limiting profiles as minimizers of a partially segregated variational problem. A preliminary study of independent interest concerns the whole-space problem. For general admissible \(K\), we determine the component-wise
Nehari infimum and prove that it is not attained, even in the radial class.
Finally, when the linear coefficients coincide, we establish a sharp criterion for the existence of positive synchronized solutions in terms of the strength of the competition.
\end{abstract}

\maketitle

\section{Introduction}

In this paper we investigate least energy fully non-trivial solutions of a
critical nonlinear Schr\"odinger system with simultaneous repulsive
interaction among all its components. We consider
\begin{equation}\label{eq:system 1} 
\begin{cases}
-\Delta u_i+\lambda_i u_i
=\mu_i|u_i|^{2^*-2}u_i
+\beta f_q(u_i)\displaystyle\prod_{j\ne i}|u_j|^q
&\text{in }\Omega\\
u_i\in H_0^1(\Omega)
\end{cases} \qquad i=1,\ldots,K,
\end{equation}
where \(\Omega\subset\mathbb R^N\) is a smooth bounded domain,
\(\mu_i>0\), \(\beta<0\), and
\[
N\ge3,\qquad K\ge3,\qquad
2^*=\frac{2N}{N-2},\qquad q=\frac{2^*}{K}\ge1.
\]
Here \(f_q(s)=|s|^{q-2}s\) for \(q>1\), while
\(f_1(s)=\operatorname{sgn}(s)\), with \(\operatorname{sgn}(0)=0\).
The identity \(Kq=2^*\) makes both the scalar nonlinearities and the
coupling Sobolev critical. The admissible pairs are
\[
\begin{array}{c|c}
N&K\\ \hline
3&3,4,5,6\\
4&3,4\\
5&3\\
6&3
\end{array}
\]
and, in particular, \(1\le q\le2\). The negative coupling penalizes simultaneous overlap of all \(K\)
components. This is a collective competition mechanism: two components
may overlap without paying an interaction cost wherever a third component
vanishes. The questions of existence and segregation must therefore be
formulated in terms of a \(K\)-fold product, and the passage from binary
to simultaneous interaction changes both the variational comparisons and
the possible limiting configurations. The cooperative case $\beta>0$ is of interest as well, but it requires completely different techniques and, for reasons of length, we have split our analysis into two papers, the present one and \cite{FaGiSo26}.

Let \(\lambda_1(\Omega)\) denote the first Dirichlet eigenvalue of
\(-\Delta\). We work in the scalar Brezis-Nirenberg range
\begin{equation}\label{eq:lambda-general}\tag{$\Lambda$}
\begin{cases}
-\lambda_1(\Omega)<\lambda_i<0,&N\ge4,\\
-\lambda_1(\Omega)<\lambda_i<-\lambda^*(\Omega),&N=3,
\end{cases}
\qquad i=1,\ldots,K,
\end{equation}
where \(\lambda^*(\Omega)\) denotes the scalar threshold in
dimension three. Thanks to \eqref{eq:lambda-general}, the quadratic forms associated with the left-hand sides in \eqref{eq:system 1} are coercive and each
uncoupled equation
\begin{equation}\label{Brezis-Nirenberg}
-\Delta v+\lambda_i v=\mu_i|v|^{2^*-2}v,
\qquad v\in H_0^1(\Omega),
\end{equation}
is a Brezis-Nirenberg problem which admits a positive ground state, whose energy is strictly below its
whole-space Sobolev level \cite{BN}. The coupled problem combines this
scalar compactness mechanism with repulsion at exactly the same critical
homogeneity.

Systems of coupled nonlinear Schr\"odinger equations are motivated,
among other contexts, by the study of interacting wave fields in
nonlinear optics. From this perspective, the components \(u_i\)
represent amplitudes of different fields, and the coefficients
\(\mu_i>0\) describe focusing self-interactions.
The product term models a collective interaction involving all
\(K\) components simultaneously, rather than a sum of pairwise
couplings. For $\beta<0$, the system can thus be viewed as an
idealized model for collective repulsive interactions, and the
choice \(Kq=2^*\) places this mechanism in the Sobolev critical regime. It is worth
mentioning that systems with multiple (i.e., beyond binary) interaction terms naturally arise in various physical contexts; see, for
instance, \cite{CCO, LQZ, Pet, Pom} and references therein.

\medskip

For two-component critical Schr\"odinger systems, Chen and Zou
\cite{Chen_Zou_2012} proved the existence of positive least energy solutions for
every negative coupling in dimension four and studied their behavior as
the competition tends to infinity. Their subsequent work \cite{Chen_Zou_2015}
treated dimensions \(N\ge5\), where the binary interaction exponent
\(2^*/2\) is smaller than two. Besides existence, it identified the
alternative between loss of one component in the weak limit and
convergence to a non-trivial segregated pair when the competition parameter $\beta \to -\infty$. In dimensions \(N\ge6\)
the first possibility is excluded, and the difference of the two limiting
components yields a least energy sign-changing solution of a
Brezis-Nirenberg-type equation. 

For systems with several components and critical pairwise coupling, we refer to
\cite{LYZ23, TY20, TYZ22} and references therein. The regularity theory for strongly
competing pairwise Schr\"odinger systems, including uniform H\"older
estimates and Lipschitz regularity of limiting profiles, was developed in
particular in \cite{NoTaTeVe, SZ15, STTZ16, ZZ26}.

Simultaneous interaction produces a different constraint. If the
competition energy is bounded along a sequence with
\(\beta\to-\infty\), a nonnegative limiting profile satisfies
\begin{equation*}
\prod_{i=1}^K u_i=0\quad\text{a.e. in }\Omega.
\end{equation*}
For \(K\ge3\), this condition does not imply
\(u_i u_j=0\) for every \(i\ne j\). It permits several components to
remain present at the same point, provided that at least one component vanishes. This distinction is reflected in the regularity and geometry of ternary and, more generally, $K$-wise competition systems and partially segregated harmonic maps studied in \cite{Gi, SoTe1,SoTe2}. These works investigate uniform estimates for singularly perturbed systems and the regularity of limiting profiles and emerging free interfaces under partial segregation. However, their results concern settings in which unconstrained energy minimizers exist. By contrast, the energy functional considered here is unbounded from below, and suitable constraints must be imposed to recover coercivity. We investigate which configurations are selected by energy minimization in the presence of Sobolev critical nonlinearities, which also introduce the possibility of concentration.

The work most closely related to the present study is that carried out by two of the authors in \cite{Giaretto_Soave_2026}, which concerns the Sobolev subcritical case. In the
competitive regime, the infimum on the component-wise Nehari set equals the sum of the scalar ground-state energies and is not
attained. Restricting to the radial class restores compactness and yields least energy fully non-trivial radial solutions. The strong-competition limit exhibits a single segregated pair of components,
while the remaining ones are positive scalar ground states, and the
difference of the pair is a least energy radial nodal solution. Our
critical problem retains the same collective interaction, but boundedness of the domain (or restriction to the radial class when working in $\R^N$) does not restore the compactness of the Sobolev embedding. Indeed, concentration phenomena create additional issues;
consequently, the existence and limiting arguments used in the subcritical regime have to be substantially rethought.

\subsection{Precise setting and main results}

As already anticipated, we are interested in \emph{fully non-trivial} solutions to \eqref{eq:system 1}, namely solutions $\mf{u}=(u_1,\dots,u_K)$ such that
$u_i\not\equiv0$ for every $i=1,\dots,K$. These have to be distinguished from \emph{semi-trivial} solutions, for which $\mf{u} \not = \boldsymbol{0}$ but $u_i \equiv 0$ for some index $i$. Note that, under assumption \eqref{eq:lambda-general}, system \eqref{eq:system 1} has several semi-trivial solutions: those obtained
by activating only the $i$-th component, 
\begin{equation}\label{semi-trivial}
\mf w_i:=(0,\dots,0,u_{\mu_i},0,\dots,0)
\end{equation}
where $u_{\mu_i}$ is a least energy positive solution to \eqref{Brezis-Nirenberg}, and also those with more than one non-trivial component, each of which is a scalar solution of the associated equation; indeed, the vanishing of a single component is enough to annihilate the full coupling term.

Our approach is variational. Set $H:=H_0^1(\Omega,\mathbb R^K)$. We define
\begin{equation}\label{energy}
	E(\mf u)
	=
	\frac12\sum_{i=1}^K
	\int_\Omega
	\left(
	|\nabla u_i|^2+\lambda_i u_i^2
	\right)\,dx
	-
	\frac{1}{2^\ast}
	\int_\Omega
	\left(
	\sum_{i=1}^K\mu_i|u_i|^{2^\ast}
	+
	K\beta\prod_{i=1}^K|u_i|^q
	\right)\,dx.
\end{equation}
A \emph{least energy fully non-trivial solution} is a fully non-trivial solution minimizing the energy among all the fully non-trivial solutions. Analogously, a
\emph{least energy positive solution} is a solution with all strictly positive components minimizing the energy among all the solutions with the same property.

A natural variational constraint to look for these solutions is the component-wise
Nehari set
\begin{equation}\label{Nehari}
	\mathcal M
	:=
	\left\{
	\mf u\in H:\ 
	u_i\not\equiv0, \ 
	\|u_i\|_{\lambda_i}^2
	=
	\mu_i\int_\Omega |u_i|^{2^\ast}\,dx
	+
	\beta\int_\Omega\prod_{j=1}^K|u_j|^q\,dx,
	\ \forall i=1,\dots,K
	\right\},
\end{equation}
where
\[
\|u_i\|_{\lambda_i}^2
:=
\int_\Omega
\left(
|\nabla u_i|^2+\lambda_i u_i^2
\right)\,dx.
\]
We define
\begin{equation}\label{eq:B definition}
	\cB
	:=
	\inf_{\mf u\in\mathcal M}E(\mf u).
\end{equation}
Summing the identities defining $\mathcal M$, we obtain that
\[
E(\mf u)
=
\frac1N
\sum_{i=1}^K
\|u_i\|_{\lambda_i}^2
\qquad
\text{for every }\mf u\in\mathcal M,
\]
and hence
\begin{equation*}
	\cB
	=
	\inf_{\mf u\in\mathcal M}
	\frac1N
	\sum_{i=1}^K
	\|u_i\|_{\lambda_i}^2.
\end{equation*}
If \(q>1\), the functional \(E\) is of class \(C^1\) on \(H\),
and its critical points are precisely the weak solutions of the system.
Moreover, the component-wise Nehari set can be equivalently written as
\[
\mathcal M
=
\left\{
\mathbf u\in H:
u_i\not\equiv 0,\;
E'(\mathbf u)[0,\ldots,0,u_i,0,\ldots,0]=0
\quad\text{for every }i=1,\ldots,K
\right\}.
\]
When \(q=1\), the interaction term is not differentiable on the whole
space \(H\). Nevertheless, the definition of \(\mathcal M\) through
the component-wise Nehari identities requires no differentiability
and remains meaningful. To obtain solutions in this case, we use
tools from nonsmooth critical point theory, which will be recalled below in Section \ref{sec:variational non smooth}.

We now state the main results of the paper. The assumption $\beta<0$ readily yields lower bounds for the norms of individual components on $\mathcal{M}$, and this allows us to prove that constrained minimizers of $E$ on $\mathcal{M}$ are actual solutions to \eqref{eq:system 1}. However, due to the critical growth, in order to prove that the infimum is attained it is necessary to rule out concentration and vanishing along constrained minimizing sequences.

We first consider the case of identical linear coefficients,
\begin{equation}\label{eq:lambda-condition}\tag{$\Lambda'$}
	\lambda_1=\cdots=\lambda_K=\lambda<0,
	\qquad
	\begin{cases}
		-\lambda_1(\Omega)<\lambda<0, & \text{if } N\geq4, \\
		-\lambda_1(\Omega)<\lambda<-\lambda^*(\Omega), & \text{if } N=3,
	\end{cases}
\end{equation}
where, as in \eqref{eq:lambda-general}, the value $\lambda^*(\Omega)\in(0,\lambda_1(\Omega))$ denotes the Brezis-Nirenberg threshold in dimension three. Let $\omega$ be a least energy positive solution of
\begin{equation}\label{BN}
	-\Delta\omega+\lambda\omega
	=
	\omega^{2^\ast-1},
	\quad
    \omega>0 \text{ in } \Omega, \quad
	\omega\in H_0^1(\Omega).
\end{equation}
A direct computation shows that $(\ell_1\omega,\dots,\ell_K\omega)$ is a positive solution of \eqref{eq:system 1} if and only if
$(\ell_1,\dots,\ell_K)\in\R_{>0}^K$ solves the algebraic system
\begin{equation}\label{eq:algebraic-system}
	\mu_i\ell_i^{2^*-2}
	+
	\beta\ell_i^{q-2}\prod_{j\neq i}\ell_j^q
	=1,
	\qquad i=1,\dots,K,
\end{equation}
where $Kq=2^*$. Solutions of this type are called \emph{synchronized solutions}. Our first result concerns the existence and uniqueness of synchronized
solutions. 
\begin{theorem}\label{thm:synchronized solutions}
	Assume \eqref{eq:lambda-condition} and let $\omega$ be a least energy positive solution of \eqref{BN}. If $\beta<0$ and
		\[
		|\beta|
		<
		\left(\prod_{i=1}^K\mu_i\right)^{1/K},
		\]
		there exists a unique
		$(\ell_1,\dots,\ell_K)\in\R_{>0}^K$ such that $(\ell_1\omega,\dots,\ell_K\omega)$
		is a positive synchronized solution of \eqref{eq:system 1}.
		If instead $\beta<0$ and
		\[
		|\beta|
		\geq
		\left(\prod_{i=1}^K\mu_i\right)^{1/K},
		\]
		then \eqref{eq:system 1} admits no positive synchronized solution.
\end{theorem}
This result is not surprising: in some sense, synchronized solutions maximize the overlap among all components and hence are not energetically favorable. The theorem shows that such solutions can exist only when the coupling parameter is sufficiently small in absolute value, namely when the competition is weak.


Now we turn to the general case where the linear coefficients
\(\lambda_i\) may differ. For each \(i=1,\ldots,K\), let
\(B_{\mu_i}\) denote the least energy level of the corresponding
scalar Brezis--Nirenberg problem \eqref{Brezis-Nirenberg}:
\[
B_{\mu_i}
:=
\frac1N
\left(
\inf_{u\in H_0^1(\Omega)\setminus\{0\}}
\frac{\|u\|_{\lambda_i}^2}
{\mu_i^{2/2^*}\|u\|_{2^*}^2}
\right)^{N/2}.
\]
Ground states of the system are nontrivial solutions of least energy,
and can be characterized as minimizers of the energy on the usual
Nehari constraint. Since \(\beta<0\), they have exactly one nonzero
component, which is a scalar ground state with energy
\(\min_{1\le j\le K}B_{\mu_j}\).
This motivates the search for least energy fully non-trivial solutions,
which we obtain by minimizing the energy on the component-wise
Nehari set. We first establish their existence in dimensions
\(N=5\) and \(N=6\).

\begin{theorem}\label{thm:least energy competitive N=5,6}
	Let $N=5$ or $N=6$, $\Omega\subset\R^N$ be a smooth bounded domain, and assume that \eqref{eq:lambda-general} holds.
	 Then, for every $\beta<0$, the level
	$\cB$ is attained by a nonnegative least energy fully non-trivial solution of
	\eqref{eq:system 1}.
\end{theorem}

In dimension $N=3$ and $N=4$, the failure of $L^2$-integrability of Talenti bubbles leads to technical complications in the proof (see Section \ref{sec:least energy competitive} below). Our argument in these dimensions requires additional assumptions, for instance on the domain, in order to recover the existence of least energy fully non-trivial solutions. We denote by $S$ the sharp Sobolev constant for the embedding
$D^{1,2}(\mathbb R^N)\hookrightarrow L^{2^\ast}(\mathbb R^N)$, namely
\begin{equation}\label{eq:sharp Sobolev R^N}
	\int_{\mathbb R^N}|\nabla u|^2\,dx
	\geq
	S
	\left(
	\int_{\mathbb R^N}|u|^{2^\ast}\,dx
	\right)^{\frac{2}{2^\ast}}
	\qquad
	\text{for every }u\in D^{1,2}(\mathbb R^N),
\end{equation}
where, as usual, $D^{1,2}(\R^N)$ represents the closure of $C^\infty_c(\R^N)$ with respect to the $L^2$-norm of the gradient. 

\begin{theorem}\label{thm:least energy competitive N=3,4}
Let $K\geq 3$, $N\in\{3,4\}$, and let
$\Omega\subset\R^N$ be a smooth bounded domain.
Assume that \eqref{eq:lambda-general} holds.
Up to relabeling the components, suppose that
\[
B_{\mu_1}\leq B_{\mu_2}\leq\cdots\leq B_{\mu_K}.
\]
Suppose in addition that
\begin{equation}\label{eq:condition-C}
\frac{1}{N}
\max_{1\leq i\leq K}\mu_i^{-\frac{N-2}{2}}
\left(
\Lambda_2(\Omega)+\max_{1\leq i\leq K}\lambda_i
\right)^{N/2}
|\Omega|
<
B_{\mu_{K-1}}
+
\frac{1}{N}
\min_{1\leq i\leq K}\mu_i^{-\frac{N-2}{2}}
S^{N/2}.
\tag{$C$}
\end{equation}
Then, for every $\beta<0$, the level
$\cB$ is reached by a nonnegative, fully non-trivial
least energy solution of \eqref{eq:system 1}.
\end{theorem}

A natural question is whether there exist domains and parameters that satisfy both \eqref{eq:lambda-general} and \eqref{eq:condition-C}. We answer this question affirmatively by exhibiting an explicit class of examples in Section \ref{sec:domains condition C}.

Concerning the asymptotic behavior as $\beta \to -\infty$, namely in the strong ``partial competition" regime, we have the following statement.

\begin{theorem}\label{thm:limit behaviour}
Suppose that $\beta_n\rightarrow-\infty$ as $n\rightarrow\infty$, and that $\mf u_n:=(u^{(n)}_1,\dots,u^{(n)}_K)\in H$ is a nonnegative least energy fully non-trivial solution of \eqref{eq:system 1} with $\beta=\beta_n$ obtained in Theorem 
\ref{thm:least energy competitive N=5,6} or in Theorem \ref{thm:least energy competitive N=3,4}. Then 
\begin{equation}\label{eq:vanishing interacion}
\displaystyle\lim_{n\rightarrow\infty}\displaystyle\int_{\Omega}\beta_n\prod_{i=1}^{K}|u^{(n)}_i|^qdx=0.
\end{equation}
Moreover, up to a subsequence if necessary, one of the following conclusions holds:
\begin{enumerate}
		\item There exists an index
        
$$
i_0 \in \argmin_{j=1,\dots,K} \left(\displaystyle\sum_{i\neq j}B_{\mu_i}+\frac{1}{N}\mu_j^{-\frac{N-2}{2}}S^{\frac N2}\right),
$$
        such that $u^{(n)}_{i_0}\rightharpoonup0$ weakly in $H_0^1(\Omega)$, and $u^{(n)}_{j}\rightarrow u^{\infty}_{j}\not=0$ strongly in
        $H_0^1(\Omega)$ for any $j\not=i_0$, where $u^{\infty}_{j}$ is a positive least energy solution of
 \begin{equation*}
 -\Delta u+\lambda_j u = \mu_j |u|^{2^\ast-2}u, \quad u\in H_0^1(\Omega);
            \end{equation*}
        \item $u^{(n)}_{i}\rightarrow u^{\infty}_{i}\not=0$ strongly in
        $H_0^1(\Omega)$ for any $i\in\{1,\cdots,K\}$, and $\mf{u}_\infty:=(u^{\infty}_{1},\dots,u^{\infty}_{K})$ is a minimizer for
    \begin{equation*}
        E_\infty(\mf{v}):= \sum_{i=1}^K \left( \frac12\int_{\Omega} (|\nabla v_i|^2+\lambda_i v_i^2)\,dx-\frac{\mu_i}{2^*}\int_{\Omega}|v_i|^{2^*}\,dx\right)
       \end{equation*}
        in the partially segregated class
        \[
        \mathcal{M}_\infty:=\left\{\mf{v} \in H\left|\begin{array}{l} v_i \not \equiv 0 \quad \forall i=1,\dots,K, \quad \prod_{i=1}^K v_i \equiv 0 \\[0.3cm]
        \|v_i\|_{\lambda_i}^2 = \mu_i \|v_i\|_{2^*}^{2^*} \quad \forall i=1,\dots,K\end{array}\right.\right\}.
        \]
        \end{enumerate}
Furthermore, if $N \in \{3,4,6\}$, only conclusion (2) holds.
\end{theorem} 

Condition \eqref{eq:condition-C} yields a uniform-in-$\beta$ energy gap below the concentration thresholds. Thus, besides ensuring existence in dimensions three and four, it prevents the loss of components in the strong competition limit. A uniform gap is also available in dimension six (on any domain), whereas in dimension five our estimates only provide a gap for each fixed coupling; this is why for $N=5$ we cannot rule out the occurrence of case (1).

\begin{remark}
  The analysis of the limiting configuration, and more generally of the partially segregated limit problem, also raises interesting open questions. On the one hand, the regularity theory developed in \cite{Gi, SoTe1, SoTe2} does not apply, since it concerns unconstrained minimizers with fixed boundary traces. On the other hand, the approach used in \cite{Giaretto_Soave_2026} cannot be applied either, since it relies crucially on the radial symmetry of the solutions.
\end{remark}

In what follows we describe the strategy of the proof of Theorems \ref{thm:least energy competitive N=5,6} and \ref{thm:least energy competitive N=3,4}, and make further comments on their differences and the role of Condition \eqref{eq:condition-C}. 
Since both the nonlinearities and the coupling term in \eqref{eq:system 1}
have critical growth, a crucial ingredient is the understanding of the whole space system
\begin{equation}\label{eq:system D12}
	\begin{cases}
		-\Delta u_i
		=
		\mu_i |u_i|^{2^\ast-2}u_i
		+
		\beta f_q(u_i)
		\displaystyle\prod_{\substack{j=1\\j\neq i}}^K |u_j|^q
		& \text{in }\mathbb R^N\\[0.4em]
		u_i\in D^{1,2}(\mathbb R^N)
	\end{cases} \qquad  i=1,\dots,K.
\end{equation}
For $\varepsilon>0$ and $y\in\mathbb R^N$, we denote by
\begin{equation}\label{attained function}
U_{\varepsilon,y}(x)
:=
[N(N-2)]^{\frac{N-2}{4}}
\left(
\frac{\varepsilon}
{\varepsilon^2+|x-y|^2}
\right)^{\frac{N-2}{2}}
\end{equation}
the Aubin--Talenti bubble \cite{Aubin_1976,Talenti_1976}. The family $\left\{
U_{\varepsilon,y}:
\varepsilon>0,\ y\in\mathbb R^N
\right\}$ exhausts all positive solutions of
\[
-\Delta u=|u|^{2^\ast-2}u
\qquad\text{in }\mathbb R^N,
\]
see \cite{CaGiSp}, each $U_{\varepsilon,y}$ attains equality in \eqref{eq:sharp Sobolev R^N}, and in particular
\begin{equation*}
	\int_{\mathbb R^N}|\nabla U_{\varepsilon,y}|^2\,dx
	=
	\int_{\mathbb R^N}|U_{\varepsilon,y}|^{2^\ast}\,dx
	=
	S^{N/2}.
\end{equation*}
As in the bounded-domain case, system \eqref{eq:system D12} also admits semi-trivial
solutions, for instance those of type
\begin{equation*}
    \mf U_i
:=
\left(
0,\dots,0,
\mu_i^{-\frac{N-2}{4}}U_{\varepsilon,y},
0,\dots,0
\right).
\end{equation*}

Set $\mathcal D
:=
D^{1,2}(\mathbb R^N,\mathbb R^K)$. We define
\begin{equation}\label{eq:energy-D12}
	I(\mf u)
	=
	\frac12
	\sum_{i=1}^K
	\int_{\mathbb R^N}|\nabla u_i|^2\,dx
	-
	\frac{1}{2^\ast}
	\int_{\mathbb R^N}
	\left(
	\sum_{i=1}^K\mu_i|u_i|^{2^\ast}
	+
	K\beta\prod_{i=1}^K|u_i|^q
	\right)\,dx.
\end{equation}
For $q>1$, $I\in C^1(\mathcal D)$. When $q=1$, we shall use
\eqref{eq:energy-D12} as the energy of the original nonnegative
problem, without claiming differentiability on the whole space.

In analogy with the problem on bounded domains, we introduce the component-wise Nehari set
\begin{equation*}
	\mathcal N
	:=
	\left\{
	\mf u\in\mathcal D:
	u_i\not\equiv0,
	\quad
	\int_{\mathbb R^N}|\nabla u_i|^2\,dx
	=
	\mu_i\int_{\mathbb R^N}|u_i|^{2^\ast}\,dx
	+
	\beta\int_{\mathbb R^N}\prod_{j=1}^K|u_j|^q\,dx,
	\ \forall i=1,\dots,K
	\right\}.
\end{equation*}
and we define
\begin{equation*}
	\cA
	:=
	\inf_{\mf u\in\mathcal N}I(\mf u).
\end{equation*}
Then we have:

\begin{theorem}\label{thm:ground-state-D12}
	Let $N\geq3$, $K\geq3$, $q=2^\ast/K\geq1$ and $\b<0$. Then,
		\begin{equation}\label{eq:A value beta negative}
		\cA
		=
		\frac1N S^{N/2}
		\sum_{i=1}^K
		\mu_i^{-\frac{N-2}{2}},
		\end{equation}
		and $\cA$ is not attained.
\end{theorem}

This identifies the whole-space energy available for a minimizing sequence whose components
all disappear in the weak limit, and plays a crucial role in the bounded-domain analysis. In fact, in proving Theorems \ref{thm:least energy competitive N=5,6} and \ref{thm:least energy competitive N=3,4}, concentration of certain components (but not of all of them) must also be excluded. Recalling from \cite{BN} that 
\[
B_{\mu_i}< \frac1N S^{N/2} \mu_i^{-\frac{N-2}{2}},
\]
the inequality which guarantees compactness of minimizing sequences turns out to be
\begin{equation}\label{eq:upper bound B}\tag{E}
	\cB<\min\left\{
		\min_{1\leq \ell\leq K}\Big(\sum_{k\not=\ell} B_{\mu_k}
		+\frac1N\mu_\ell^{-\frac{N-2}{2}}S^{\frac N2}\Big),\cA
		\right\}.
\end{equation}
The mixed levels represent configurations consisting of $K-1$ surviving
scalar components and a concentrating $K$-th component. A bound below
the value $\mathcal A$ alone would not rule out this loss of compactness. We prove
that \eqref{eq:upper bound B} yields attainment of
$\mathcal{B}$, using splitting of both the scalar critical terms and the
product interaction, in the spirit of the Brezis--Lieb lemma
\cite{BL83}. Although an estimate similar to \eqref{eq:upper bound B} was also used to treat the case of the binary interaction, verifying its validity requires entirely new arguments; 
The main difference from the binary case lies in the estimates of the interaction term, which require competitors and concentration scales adapted to the simultaneous coupling.

In dimensions \(N=5,6\), we establish these strict inequalities for every negative coupling, independently of the domain. We construct test functions by combining two scalar ground states with a truncated critical bubble placed near the boundary. Both scalar factors are then small near the bubble, allowing us to control the simultaneous interaction. The bubble scale and its distance from the boundary are chosen so that the interaction and truncation errors are dominated by the gain arising from the negative linear coefficient. 

For \(N=3,4\), we give sufficient conditions for the same compactness inequality. The previous argument fails because of the lack of $L^2$-integrability of the Aubin--Talenti bubble in $\R^N$. However, a useful class of competitors consists of those having two functions
with disjoint supports, and the remaining $K-2$ scalar ground states: indeed, separation of
just a single pair eliminates the entire interaction, independently of
\(\beta\). In order to make some quantitative estimates, we choose as a segregated pair the positive and negative parts of a second Dirichlet
eigenfunction; this gives the sufficient condition \eqref{eq:condition-C} in Theorem \ref{thm:least energy competitive N=3,4}. It is worth mentioning that this condition is certainly not sharp. Other possible conditions can be obtained by choosing different segregated pairs, such as the positive and negative parts of a least energy nodal solution to \eqref{Brezis-Nirenberg}. However, such computations become more involved and it is not clear that there is an optimal choice. Consequently, we present the simpler estimate \eqref{eq:condition-C}, which has the merit of involving only the data in an explicit, and hence easily verifiable, way; as already mentioned, we refer to Section \ref{sec:domains condition C} for more details.

\medskip

Finally, returning to the whole-space problem \eqref{eq:system D12}, since the least energy level $\cA$ is not attained, a natural question is whether there exists a least energy fully non-trivial radial solution for \eqref{eq:system D12}. This is the case in the subcritical regime \cite{Giaretto_Soave_2026}, where one can minimize in the radial component-wise Nehari set. Here, due to possible concentration phenomena, such a minimization is instead not fruitful: precisely, defining $\mathcal N_{\rad} = \mathcal{N} \cap D^{1,2}_{\rad}(\R^N)$ and $\cA_{\rad}$ the associated infimum level, the following holds:

\begin{theorem}\label{thm:rad-D12}
Let $N\geq3$, $K\geq3$, $q=2^\ast/K\geq1$ and $\b<0$. Then,
		\begin{equation*}
		\cA_{\rad}=\cA,
		\end{equation*}
		and $\cA_{\rad}$ is not attained.
\end{theorem}

Regarding the existence of positive solutions, we have a precise counterpart of Theorem \ref{thm:synchronized solutions} in the whole space setting, ensuring existence for sufficiently small coupling parameters.

\begin{theorem}\label{thm:rad-synchronized solutions}
	Let $U=U_{\eps,y}$ be an Aubin--Talenti bubble. If $\beta<0$ and
		\[
		|\beta|
		<
		\left(\prod_{i=1}^K\mu_i\right)^{1/K},
		\]
		there exists a unique
		$(\ell_1,\dots,\ell_K)\in\R_{>0}^K$ such that	$(\ell_1U,\dots,\ell_KU)$
		is a positive synchronized solution of \eqref{eq:system D12}.
		If instead  $\beta<0$ and
		\[
		|\beta|
		\geq
		\left(\prod_{i=1}^K\mu_i\right)^{1/K},
		\]
		then \eqref{eq:system D12} admits no positive synchronized solution.
\end{theorem}

\subsection*{Structure of the paper}
We conclude the introduction with a brief description of the organization
of the paper. In Section~\ref{sec:synchronized}, we prove the results concerning synchronized solutions, namely Theorems \ref{thm:synchronized solutions} and \ref{thm:rad-synchronized solutions}. In Section~\ref{sec:auxiliary},
we analyze the limiting critical problem in $\mathbb R^N$.
In Section~\ref{sec:variational non smooth}, we recall the tools from nonsmooth critical point theory needed to treat the case $q=1$, while in Section~\ref{sec:least energy competitive} we prove Theorems~\ref{thm:least energy competitive N=5,6} and \ref{thm:least energy competitive N=3,4}. In Section \ref{sec:domains condition C} we provide examples of domains satisfying condition \eqref{eq:condition-C}, while in Section \ref{sec:asymptotic behaviour} we prove Theorem \ref{thm:limit behaviour}.

\section{Synchronized solutions}\label{sec:synchronized}
The present section is devoted to the proof of Theorem \ref{thm:synchronized solutions} and Theorem \ref{thm:rad-synchronized solutions}.
\begin{proof}[Proof of existence and uniqueness part in Theorem \ref{thm:synchronized solutions}.]
As observed in the introduction, it is enough to prove the corresponding statement
	for the algebraic system \eqref{eq:algebraic-system}. For
	$\tau\in[0,1]$, consider the equation
	\begin{equation}\label{eq:algebraic-homotopy}
		G_i(\ell_1,\dots,\ell_K)
		:=
		\mu_i\ell_i^{2^*}
		+
		\tau\beta\prod_{j=1}^K\ell_j^q
		-
		\ell_i^2
		=
		0,
		\qquad i=1,\dots,K,
	\end{equation}
	and set
	\[
	Z
	:=
	\left\{
	\tau\in[0,1]:
	\eqref{eq:algebraic-homotopy}
	\text{ has a unique solution in }\R_{>0}^K
	\right\}.
	\]
	We prove that $Z=[0,1]$.
	
	\medskip
	\noindent\textit{Step 1: $0\in Z$.}
	For $\tau=0$, the system is decoupled and reduces to
	\[
	\mu_i\ell_i^{2^*}-\ell_i^2=0,
	\qquad i=1,\dots,K.
	\]
	Hence its unique positive solution is
	\[
	\overline\ell_i
	=
	\mu_i^{-\frac{1}{2^*-2}},
	\qquad i=1,\dots,K.
	\]
	Therefore $0\in Z$.
	
	\medskip
	\noindent\textit{Step 2: $Z$ is open in $[0,1]$.}
	Let $\tau_0\in Z$ and let
	$(\ell_1^0,\dots,\ell_K^0)\in\R_{>0}^K$ be the unique solution of
	\eqref{eq:algebraic-homotopy} for $\tau=\tau_0$. Consider the
	Jacobian matrix
	\[
	M
	= (M_{ij})_{i,j=1,\dots,K}=
	\left(
	\partial_{\ell_j}G_i
	\right)_{i,j=1,\dots,K}.
	\]
	Using \eqref{eq:algebraic-homotopy}, we obtain
	\[
	M_{ii}
	=
	(\ell_i^0)^{-1}
	\left[
	(2^*-2)\mu_i(\ell_i^0)^{2^*}
	+
	(q-2)\tau_0\beta
	\prod_{j=1}^K(\ell_j^0)^q
	\right]
	\]
	and, for $i\neq j$,
	\[
	M_{ij}
	=
	(\ell_j^0)^{-1}
	\tau_0q\beta
	\prod_{h=1}^K(\ell_h^0)^q.
	\]
	Thus
	\[
	\det M
	=
	\prod_{i=1}^K(\ell_i^0)^{-1}
	\det\widetilde M,
	\]
	where
	\[
	\widetilde M_{ii}
	=
	(2^*-2)\mu_i(\ell_i^0)^{2^*}
	+
	(q-2)\tau_0\beta
	\prod_{j=1}^K(\ell_j^0)^q
	\]
	and
	\[
	\widetilde M_{ij}
	=
	\tau_0q\beta
	\prod_{j=1}^K(\ell_j^0)^q,
	\qquad i\neq j.
	\]
	We claim that
	\begin{equation}\label{eq:jacobian-nondegenerate}
		\det M\neq0.
	\end{equation}
Indeed, since $\beta<0$ and
	$G_i(\ell_1^0,\dots,\ell_K^0)=0$, we have
	\[
	\mu_i(\ell_i^0)^{2^*}
	=
	(\ell_i^0)^2
	+
	\tau_0|\beta|
	\prod_{j=1}^K(\ell_j^0)^q
	>
	\tau_0|\beta|
	\prod_{j=1}^K(\ell_j^0)^q.
	\]
	Hence
	\begin{equation*}
		|\tilde M_{ii}|=(2^*-2)\mu_i(\ell_i^0)^{2^*}
		-
		(q-2)\tau_0|\beta|
		\prod_{j=1}^K(\ell_j^0)^q
		>
		(K-1)q\tau_0|\beta|
		\prod_{j=1}^K(\ell_j^0)^q=\sum_{j\neq i} |\tilde M_{ij}|,
	\end{equation*}
	where we used $Kq=2^*$. Thus $\widetilde M$ is strictly
	diagonally dominant and therefore nonsingular. This proves
	\eqref{eq:jacobian-nondegenerate}.

By \eqref{eq:jacobian-nondegenerate}, the Implicit Function Theorem
yields an open neighborhood $U_0$ of $\tau_0$ and a neighborhood
$A_0\subset\R_{>0}^K$ of
$(\ell_1^0,\dots,\ell_K^0)$ such that
\eqref{eq:algebraic-homotopy} is uniquely solvable in
$U_0\times A_0$.

Suppose, by contradiction, that there exists $\tau_1\in U_0$ such that
\eqref{eq:algebraic-homotopy} admits a second solution
\[
(\widehat\ell_1^0,\dots,\widehat\ell_K^0)
\in
\R_{>0}^K\setminus A_0.
\]
By an analogous reasoning applying the Implicit Function Theorem,
there exist $\varepsilon>0$ and a solution curve
\[
\tau
\longmapsto
\bigl(
\widehat\ell_1(\tau),\dots,\widehat\ell_K(\tau)
\bigr)
\]
defined for
$\tau\in(\tau_1-\varepsilon,\tau_1+\varepsilon)$, satisfying
\eqref{eq:algebraic-homotopy} and such that
\[
\bigl(
\widehat\ell_1(\tau_1),\dots,\widehat\ell_K(\tau_1)
\bigr)
=
(\widehat\ell_1^0,\dots,\widehat\ell_K^0).
\]
Without loss of generality, assume that $\tau_0<\tau_1$, and extend this
curve to the left as long as possible. Since the curve cannot reach
$\tau_0$ while entering $U_0\times A_0$, there exists
$\tau_2\in[\tau_0,\tau_1)$ such that
\[
\bigl(
\widehat\ell_1(\tau),\dots,\widehat\ell_K(\tau)
\bigr)
\]
is defined for $\tau\in(\tau_2,\tau_1]$, while, as
$\tau\to\tau_2^+$, either at least one component becomes unbounded or
the curve approaches the boundary of $\R_{>0}^K$.

Noticing $\beta<0$ and
\[
|\beta|
<
\left(
\prod_{i=1}^K\mu_i
\right)^{1/K}.
\]
We claim that there exists $D_0>0$ such that
\begin{equation}\label{eq:uniform-upper-negative-beta}
\widehat\ell_i(\tau)<D_0,
\qquad i=1,\dots,K,
\end{equation}
as $\tau\to\tau_2^+$. Suppose otherwise. Then, up to a subsequence,
there exists $i_0\in\{1,\dots,K\}$ such that
\[
\widehat\ell_{i_0}(\tau)\to+\infty
\qquad\text{as }\tau\to\tau_2^+,
\]
and we may choose $i_0$ so that
\[
\frac{\widehat\ell_j(\tau)}
{\widehat\ell_{i_0}(\tau)}
\leq C
\]
for every $j=1,\dots,K$ and $\tau$ sufficiently close to $\tau_2$.

Assume first that, for some
$j_0\in\{1,\dots,K\}\setminus\{i_0\}$,
\[
\frac{\widehat\ell_{j_0}(\tau)}
{\widehat\ell_{i_0}(\tau)}
\to0
\qquad\text{as }\tau\to\tau_2^+.
\]
Dividing the $i_0$-th equation in
\eqref{eq:algebraic-homotopy} by
$\widehat\ell_{i_0}(\tau)^{2^*}$ and recalling that $Kq=2^*$, we obtain
\[
\mu_{i_0}
=
-\frac{\tau\beta}{\widehat\ell_{i_0}(\tau)^{2^*}}\prod_{j=1}^K\widehat\ell_j(\tau)^q
+
\widehat\ell_{i_0}(\tau)^{2-2^*}
\longrightarrow0,
\]
as $\tau\to\tau_2^+$, which is impossible.

Consequently, up to a subsequence, there exist $c_j>0$ such that
\[
\frac{\widehat\ell_j(\tau)}
{\widehat\ell_{i_0}(\tau)}
\longrightarrow c_j>0,
\qquad j=1,\dots,K.
\]
In particular, all components tend to $+\infty$.
Dividing the $i$-th equation in
\eqref{eq:algebraic-homotopy} by
$\widehat\ell_i(\tau)^{2^*}$, we obtain
\[
\mu_i
=
-\tau\beta
\frac{\prod_{j=1}^K\widehat\ell_j(\tau)^q}
{\widehat\ell_i(\tau)^{2^*}}
+
o_\tau(1),
\qquad i=1,\dots,K,
\]
where $o_\tau(1)\to 0$ as $\tau\to\tau_2^+$.
Multiplying these identities over $i=1,\dots,K$, using
$Kq=2^*$ and then letting $\tau\to\tau_2^+$ gives
\[
\prod_{i=1}^K\mu_i
=
|\tau_2\beta|^K
\leq
|\beta|^K,
\]
which contradicts
\[
|\beta|
<
\left(
\prod_{i=1}^K\mu_i
\right)^{1/K}.
\]
This proves \eqref{eq:uniform-upper-negative-beta}.

On the other hand, from \eqref{eq:algebraic-homotopy} we have, when
$\beta<0$,
\[
\widehat\ell_i(\tau)^2
\leq
\mu_i\widehat\ell_i(\tau)^{2^*},
\]
and therefore
\[
\widehat\ell_i(\tau)
\geq
\mu_i^{-\frac{1}{2^*-2}},
\qquad i=1,\dots,K.
\]
Thus also in this case the solution curve remains in a compact subset
of $\R_{>0}^K$, which contradicts its maximality.

We conclude that no second positive solution can occur for
$\tau\in U_0$. Hence $U_0\subset Z$, and $Z$ is open in $[0,1]$.

\medskip

\noindent\textit{Step 3: $Z$ is closed in $[0,1]$.}
Let $\{\tau_n\}\subset Z$ be such that
$\tau_n\to\tau_0\in[0,1]$, and let $(\ell_1^n,\dots,\ell_K^n)\in\R_{>0}^K$
be the unique positive solution of
\eqref{eq:algebraic-homotopy} corresponding to $\tau=\tau_n$.
By the estimates established above, the sequence
$\{(\ell_1^n,\dots,\ell_K^n)\}$ is bounded. Hence, up to a subsequence,
\[
(\ell_1^n,\dots,\ell_K^n)
\longrightarrow
(\ell_1^0,\dots,\ell_K^0)
\in\R_+^K.
\]
Passing to the limit in \eqref{eq:algebraic-homotopy}, we see that
$(\ell_1^0,\dots,\ell_K^0)$ solves the same system with
$\tau=\tau_0$. We are left to prove that such solution is unique.

Actually, since $\beta<0$, it follows directly from
\eqref{eq:algebraic-homotopy} that
\[
(\ell_i^n)^2
\leq
\mu_i(\ell_i^n)^{2^*},
\]
and hence
\[
\ell_i^0
\geq
\mu_i^{-\frac{1}{2^*-2}}
>0,
\qquad i=1,\dots,K.
\]
Thus, in either case, $(\ell_1^0,\dots,\ell_K^0)\in\R_{>0}^K$. Arguing as before with the Implicit Function Theorem,
the limiting solution is locally unique. If another positive solution
existed at $\tau=\tau_0$, the same continuation argument used in
Step~2 would produce a second positive branch for $\tau_n$ sufficiently
close to $\tau_0$, contradicting $\tau_n\in Z$. Hence
$(\ell_1^0,\dots,\ell_K^0)$ is the unique positive solution at
$\tau=\tau_0$, and therefore $\tau_0\in Z$.

This proves that $Z$ is closed in $[0,1]$. Since $Z$ is nonempty, open,
and closed in the connected set $[0,1]$, we conclude that $Z=[0,1]$. In particular, taking $\tau=1$, the algebraic system
\eqref{eq:algebraic-system} has a unique solution in $\R_{>0}^K$.
\end{proof}

It remains to prove the nonexistence assertion in Theorem \ref{thm:synchronized solutions}.
\begin{proof}[Proof of nonexistence part in Theorem \ref{thm:synchronized solutions}]
Assume that
\[
|\beta|
\geq
\left(
\prod_{i=1}^K\mu_i
\right)^{1/K},
\]
and suppose, by contradiction, that
$(\ell_1^0,\dots,\ell_K^0)\in\R_{>0}^K$ solves
\eqref{eq:algebraic-system}. Then
\[
(\ell_i^0)^2
+
|\beta|
\prod_{j=1}^K(\ell_j^0)^q
=
\mu_i(\ell_i^0)^{2^*},
\qquad i=1,\dots,K.
\]
In particular,
\[
|\beta|
\prod_{j=1}^K(\ell_j^0)^q
<
\mu_i(\ell_i^0)^{2^*},
\qquad i=1,\dots,K.
\]
Multiplying these inequalities over $i=1,\dots,K$ and using
$Kq=2^*$, we obtain
\[
|\beta|^K
\prod_{j=1}^K(\ell_j^0)^{2^*}
<
\left(
\prod_{i=1}^K\mu_i
\right)
\prod_{j=1}^K(\ell_j^0)^{2^*}.
\]
Since all $\ell_j^0$ are positive, this yields
\[
|\beta|
<
\left(
\prod_{i=1}^K\mu_i
\right)^{1/K},
\]
which contradicts the assumption. This completes the proof.
\end{proof}

Finally, we give the proof of Theorem \ref{thm:rad-synchronized solutions}.
\begin{proof}[Proof of Theorem \ref{thm:rad-synchronized solutions}]
Let $U$ be an Aubin--Talenti bubble. Then $(\ell_1U,\dots,\ell_KU)$ is a positive synchronized solution of \eqref{eq:system D12} if and only if
$(\ell_1,\dots,\ell_K)\in\R_{>0}^K$ solves the algebraic system
\begin{equation*}
	\mu_i\ell_i^{2^*-2}
	+
	\beta\ell_i^{q-2}\prod_{j\neq i}\ell_j^q
	=1,
	\qquad i=1,\dots,K.
\end{equation*}
The remainder of the proof is therefore the same as that of Theorem~\ref{thm:synchronized solutions}, and so we omit it here.
\end{proof}

\section{The limiting problem in \texorpdfstring{$\R^N$}{RN}}\label{sec:auxiliary}
The present section is devoted to the proof of the results concerning the limiting problem in the whole space, collected in Theorem \ref{thm:ground-state-D12} and \ref{thm:rad-D12}.

\begin{proof}[Proof of Theorem \ref{thm:ground-state-D12}]
	Define
	\[
	(u_{1,R}(x),\ldots,u_{K,R}(x))
	=
	\bigl(w_{\mu_1}(x+Re_1),w_{\mu_2}(x),\ldots,w_{\mu_K}(x)\bigr),
	\]
	where $w_{\mu_i}$ is a positive least energy solution of
	\[
	-\Delta u=\mu_i u^{2^*-1}
	\qquad\text{in }D^{1,2}(\R^N).
	\]
	Then
	\[
	\lim_{R\to+\infty}
	\int_{\R^N}\prod_{i=1}^K u_{i,R}^{q}dx
	=0.
	\]
	
	We prove that, for $R$ large enough, there exist
	$t_{i,R}>0$, $i=1,\ldots,K$, such that
	\[
	F_{i,R}(t_1,\dots, t_K)
	:=
	t_{i}^2\int_{\R^N}|\nabla u_{i,R}|^2dx
	-
	t_{i}^{2^*}\mu_i\int_{\R^N}|u_{i,R}|^{2^*}dx
	-
	t_{i}^{q}
	\beta
	\prod_{j\not=i} t_{j}^{q}
	\int_{\R^N}\prod_{j=1}^K u_{j,R}^{q}dx
	=0.
	\]
	
	Fix $\ell>0$. For $t_i=\ell$, we have, since $\b<0$,
	\[
	F_{i,R}\big|_{t_i=\ell}
	=
	\ell^2\int_{\R^N}|\nabla u_{i,R}|^2dx
	-
	\ell^{2^*}\mu_i\int_{\R^N}|u_{i,R}|^{2^*}dx
	-
	\ell^{q}
	\beta \prod_{j\neq i}t_j^{q}
	\int_{\R^N}\prod_{j=1}^K u_{j,R}^{q}dx
	\geq
	C_{i}\ell^2 -C_{i} \ell^{2^*},
	\]
	with $C_{i}=\norm{\nabla w_{\mu_i}}_2^2=\mu_i \norm{w_{\mu_i}}_{2^*}^{2^*}$.
	Therefore we can find $\bar\ell>0$ small enough such that, for every $i$,
	$$
	F_{i,R}\big|_{t_i=\bar\ell}>0.
	$$
	
	On the other hand, fix $L>\bar\ell$, then for $t_i=L$, $t_j\in [\bar\ell,L]$ for every $j\neq i$, 
	\[
	\begin{aligned}
	F_{i,R}\big|_{t_i=L}
	&\leq
	L^2\int_{\R^N}|\nabla u_{i,R}|^2dx
	-
	L^{2^*}\mu_i\int_{\R^N}|u_{i,R}|^{2^*}dx
	+
	L^{q}|\beta|
	\prod_{j\neq i}t_j^{q}
	\int_{\R^N}\prod_{j=1}^K u_{j,R}^{q}dx \\
	&\leq L^2\int_{\R^N}|\nabla u_{i,R}|^2dx
	-
	L^{2^*}\mu_i\int_{\R^N}|u_{i,R}|^{2^*}dx
	+
	L^{2^*}|\beta|
	\int_{\R^N}\prod_{j=1}^K u_{j,R}^{q}dx 
	\end{aligned}
	\]
	and hence
	\[
	F_{i,R}\big|_{t_i=L}
	\leq
	L^2 C_i-L^{2^*}(C_i+o(1)),
	\]
	with $o(1)\to0$ as $R\to+\infty$. Therefore, for every $R>\bar R$ for some $\bar R$ large enough, we can find $\bar L$ (independent of $R$) such that
	\[
	F_{i,R}\big|_{t_i=\bar L}\leq C_i \bar L^2-\frac{C_i}{2}\bar L^{2^*}<0.
	\]
		
	Thus, by Miranda's theorem, for every $R>\bar R$ there exist
	$t_{i,R}\in[\bar \ell,\bar L]$, $i=1,\ldots,K$, such that
	\[
	F_{i,R}(t_{1,R},\ldots,t_{K,R})=0,
	\qquad i=1,\ldots,K,
	\]
	namely 
	$$
	(t_{1,R}u_{1,R},\dots, t_{K,R}u_{K,R})\in \cN.
	$$
	
	Moreover,
	\[
	t_{i,R}^2
	\int_{\R^N}|\nabla w_{\mu_i}|^2
	=
	t_{i,R}^{2^*}
	\mu_i\int_{\R^N}|w_{\mu_i}|^{2^*}
	+
	\beta \,
	t_{i,R}^{q}
	\prod_{j\not=i}^K t_{j,R}^{q}
	\int_{\R^N}
	\prod_{j=1}^K u_{j,R}^{q} 
	\]
	Since the last integral tends to zero as $R\to+\infty$ and $t_{i,R}$ is bounded, we obtain
	\[
	t_{i,R}\to1
	\qquad\text{as }R\to+\infty.
	\]
	
	Therefore,
	\[
	\cA
	\leq
	I(t_{1,R}u_{1,R},\ldots,t_{K,R}u_{K,R})
	=
	\frac1N
	\sum_{i=1}^K
	t_{i,R}^2
	\int_{\R^N}|\nabla u_{i,R}|^2.
	\]
	Passing to the limit as $R\to+\infty$, we obtain
	\[
	\cA
	\leq
	\frac1N
	S^{N/2}
	\sum_{i=1}^K
	\mu_i^{-\frac{N-2}{2}}.
	\]
	
	On the other hand, for any $\mf u=(u_1,\dots,u_K)\in\mathcal N$, we see from $\beta<0$ and \eqref{eq:sharp Sobolev R^N} that, for every $i=1,\dots,K$,
	\[
	\int_{\R^N} |\nabla u_i|^2\,dx
	\leq
	\mu_i\int_{\R^N}|u_i|^{2^*}\,dx
	\leq
	\mu_i S^{-2^*/2}
	\left(
	\int_{\R^N}|\nabla u_i|^2\,dx
	\right)^{2^*/2},
	\]
	and hence
	\beq\label{10111}
	\int_{\R^N}|\nabla u_i|^2\,dx
	\geq
	\mu_i^{-\frac{N-2}{2}}S^{N/2}.
	\eeq
	Combining these inequalities with the expression of the functional on $\mathcal{N}$, we obtain that
	\[
	\cA
	\geq
	\frac{1}{N}
	\sum_{i=1}^K
	\mu_i^{-\frac{N-2}{2}}S^{N/2},
	\]
	whence the equality follows.
    
    Now, assume by contradiction that $\cA$ is attained by some
	$\mf u=(u_1,\dots,u_K)\in\mathcal N$. Since
	\[
	\cA
	=
	\frac1N\sum_{i=1}^K\int_{\R^N}|\nabla u_i|^2\,dx
	=
	\frac1N
	\left(
	\sum_{i=1}^K\mu_i^{-\frac{N-2}{2}}
	\right)
	S^{N/2},
	\]
	and \eqref{10111} holds for each $u_i$,	we necessarily have
	\[
	\int_{\R^N}|\nabla u_i|^2\,dx
	=
	\mu_i^{-\frac{N-2}{2}}S^{N/2}
	\qquad\text{for every }i=1,\dots,K.
	\]
	Hence equality must hold in
	\[
	\int_{\R^N}|\nabla u_i|^2\,dx
	\leq
	\mu_i\int_{\R^N}|u_i|^{2^*}\,dx
	\leq
	\mu_i S^{-2^*/2}
	\left(
	\int_{\R^N}|\nabla u_i|^2\,dx
	\right)^{2^*/2}.
	\]
	This, together with the
	Nehari identity, yields
	\[
	\int_{\R^N}\prod_{j=1}^K |u_j|^{q}\,dx=0.
	\]
	On the other hand, equality in the Sobolev inequality implies that
	each $u_i$ is an Aubin--Talenti instanton.
	Therefore $|u_i|>0$ in $\R^N$ for every $i=1,\dots,K$, and thus the interaction term must be strictly positive, which is a contradiction.
\end{proof}
\begin{proof}[Proof of Theorem \ref{thm:rad-D12}]
For any $\epsilon>0$, define $U_{\epsilon,0}$ as in \eqref{attained function} with $y=0$.
Set $a_i:=\mu_i^{-\frac{N-2}{4}}$, $i=1,\cdots,K$. Then
\begin{equation*}
 -\Delta (a_iU_{\varepsilon,0})
=
\mu_i (a_iU_{\varepsilon,0})^{2^\ast-1}
\qquad\text{in }\mathbb R^N   
\end{equation*}
and
\begin{equation}\label{eq:2}
	\int_{\mathbb R^N}|\nabla (a_iU_{\varepsilon,0})|^2\,dx
	=
	\mu_i \int_{\mathbb R^N}|(a_iU_{\varepsilon,0})|^{2^\ast}\,dx
	=
	\mu_i^{-\frac{N-2}{2}}S^{N/2}.
\end{equation}
Set $\mf v^{\epsilon}:=(v_1^{\epsilon},\cdots,v_K^{\epsilon})$, where $v_1^{\epsilon}:=a_1U_{\varepsilon,0}$, and $v_i^{\epsilon}:=a_iU_{1,0}$ for $i=2,\cdots,K$.
From \eqref{eq:2}, we have
\begin{equation}\label{eq:3}
	\int_{\mathbb R^N}|\nabla v_i^{\epsilon}|^2\,dx
	=
	\mu_i\int_{\mathbb R^N}|v_i^{\epsilon}|^{2^\ast}\,dx
	=
	\mu_i^{-\frac{N-2}{2}}S^{N/2},
\end{equation}
for $i=1,\cdots,K$. We claim that 
\begin{equation}\label{eq:4}
	K_\epsilon:=
	\int_{\mathbb R^N}\displaystyle\prod_{i=1}^K|v_i^{\epsilon}|^qdx\rightarrow0 \qquad \text{for }\epsilon\rightarrow0^+.
\end{equation}
Indeed, up to a positive constant independent of $\epsilon$,
$$K_\epsilon\leq C\displaystyle\int_{0}^{\infty}\left(\frac{\epsilon}{\epsilon^2+r^2}\right)^{d_1}\prod_{i=2}^{K}\left(\frac{1}{1+r^2}\right)^{d_i}r^{N-1}dr,$$
where $d_i:=\frac{N}{K}$ for $i=1,\cdots,K$. For $r\geq1$,
\begin{equation*}
\displaystyle\int_{1}^{\infty}\left(\frac{\epsilon}{\epsilon^2+r^2}\right)^{\frac{N}{K}}\prod_{i=2}^{K}\left(\frac{1}{1+r^2}\right)^{\frac{N}{K}}r^{N-1}dr\leq C\epsilon^{\frac{N}{K}}\displaystyle\int_{1}^{\infty} r^{-N-1}dr\rightarrow0,
\end{equation*}
as $\epsilon\rightarrow 0^+$. 
For $0<r\leq1$, the product over $i=2,\cdots,K$ is bounded, and
\begin{align*}
\displaystyle\int_{0}^{1}\left(\frac{\epsilon}{\epsilon^2+r^2}\right)^{\frac{N}{K}}\prod_{i=2}^{K}\left(\frac{1}{1+r^2}\right)^{\frac{N}{K}}r^{N-1}dr &\leq C\displaystyle\int_{0}^{1}\left(\frac{\epsilon}{\epsilon^2+r^2}\right)^{\frac{N}{K}}r^{N-1}dr
\\&=C\epsilon^{\frac{K-1}{K}N}\int_{0}^{\frac{1}{\epsilon}}\frac{t^{N-1}}{(1+t^2)^\frac{N}{K}}dt \approx C\epsilon^{\frac{K-1}{K}N}+\epsilon^\frac{N}{K}\rightarrow0,
\end{align*}
as $\epsilon\rightarrow0^+$. Hence, \eqref{eq:4} holds.

Now we project $\mf v^{\epsilon}$ onto  $\mathcal N_{\rad}$.
For $\mf t=(t_1,\cdots,t_K)\in\mathbb{R}_+^K$ and $\kappa\in \mathbb{R}$, define
\begin{equation*}
\Phi_i(\mf t,\kappa):=t_i^2 A_i-t_i^{2^\ast} A_i-\beta\displaystyle\prod_{j=1}^{K}t_j^q \kappa, ~i=1,\cdots,K,
\end{equation*}
where
\begin{equation*}
 A_i=
	\int_{\mathbb R^N}|\nabla v_i^{\epsilon}|^2\,dx
	=
	\mu_i\int_{\mathbb R^N}|v_i^{\epsilon}|^{2^\ast}\,dx
	=
	\mu_i^{-\frac{N-2}{2}}S^{N/2},
\end{equation*}
by \eqref{eq:3}. At $\kappa=0$,
\begin{equation*}
\Phi_i(\mf t,0):= A_it_i^2(1-t_i^{2^\ast-2}), 
\end{equation*}
so the unique positive zero is $t_i=1$ for all $i\in\{1,\cdots,K\}$. Moreover,
\begin{equation*}
 \left(\frac{\partial\Phi_i}{\partial t_j}\right)_{i,j}=\operatorname{diag}((2-2^*)A_1,\cdots,(2-2^*)A_K)
\end{equation*}
at $\mf t=(1,\cdots,1)$, which is invertible since $2^*>2$. Applying the Implicit function theorem, there is $\delta>0$ and a unique $C^1$ map
\begin{equation*}
 \mf t(\kappa)=(t_1(\kappa),\cdots,t_K(\kappa)),~|\kappa|<\delta,
\end{equation*}
such that
\begin{equation*}
 \mf t(0)=(1,\cdots,1),~\Phi(\mf t(\kappa),\kappa)=0.
\end{equation*}
For $\epsilon>0$ small enough, we have $|K_\epsilon|<\delta$. Set
\begin{equation*}
t_i^{\epsilon}:=t_i(K_\epsilon),~i=1,\cdots,K.
\end{equation*}
Then
\begin{equation*}
t_i^{\epsilon}\rightarrow1,~\epsilon\rightarrow0^+,
\end{equation*}
and
\begin{equation*}
(t_i^{\epsilon})^2A_i=(t_i^{\epsilon})^{2^*}A_i+\beta\displaystyle\displaystyle\prod_{j=1}^{K}(t^{\epsilon}_j)^qK_\epsilon.
\end{equation*}
for $i=1,\cdots, K$. Therefore, $\mf u_\epsilon:=(t_1^{\epsilon}v_1^{\epsilon},\cdots,t_K^{\epsilon}v_K^{\epsilon})\in \mathcal N_{rad}$, and so
\begin{align*}
I(\mf u_\epsilon)&=\frac{1}{N}\displaystyle\sum_{i=1}^{K}\int_{\mathbb R^N}|\nabla u_i^{\epsilon}|^2\,dx=\frac{1}{N}\displaystyle\sum_{i=1}^{K}(t_i^{\epsilon})^2A_i
=\frac{1}{N}\displaystyle\sum_{i=1}^{K}(t_i^{\epsilon})^2\mu_i^{-\frac{N-2}{2}}S^{N/2}
\\&\rightarrow \frac{1}{N}\displaystyle\sum_{i=1}^{K}\mu_i^{-\frac{N-2}{2}}S^{N/2}=\cA,
\end{align*}
as $\epsilon\rightarrow0^+$. Thus,
\begin{equation*}
\cA_{\rad}\leq\cA.
\end{equation*}
On the other hand, $\cA_{\rad}\geq\cA$ since $\mathcal N_{\rad}\subset \mathcal N$, we obtain that
\eqref{eq:A value beta negative} holds. Furthermore, since $\cA$ cannot be attained
by Theorem \ref{thm:ground-state-D12}, we obtain that $\cA_{\rad}$ cannot be attained either. This completes the proof.
\end{proof}

\section{Variational setting in the non-differentiable case \texorpdfstring{$q=1$}{q=1}}\label{sec:variational non smooth}

In the borderline case $q=1$, the functional $E$ is only locally
Lipschitz on $H$. We briefly recall the nonsmooth variational framework
that will be used below. We refer to \cite{Chang_1981,Clarke_1976} for the
general theory and to \cite{Clapp_Saldaña_Schiera_2026} for a closely
related argument.

Let $X$ be a Banach space and let $F:X\to\R$ be locally Lipschitz.
The generalized gradient of $F$ at $u\in X$ is defined by
\[
\partial F(u)
:=
\left\{
g\in X^*:
\langle g,v\rangle
\leq
\limsup_{\substack{w\to u\\ t\to0^+}}
\frac{F(w+tv)-F(w)}{t}
\quad\text{for every }v\in X
\right\},
\]
where $\langle\cdot,\cdot\rangle$ denotes the usual duality pairing.
We recall that
\[
\partial F(u)=\{F'(u)\}
\quad\text{if }F\in\mathcal C^1(X),
\]
and that
\[
\partial(F_1+F_2)(u)
\subseteq
\partial F_1(u)+\partial F_2(u),
\qquad
\partial(\lambda F)(u)
=
\lambda\partial F(u).
\]

We now specialize these notions to the functional $E$ in
\eqref{energy}. Throughout this section we assume $q=1$, and hence
$K=2^*$. For $t\in\R$, set
\[
\operatorname{Sign}(t)
:=
\begin{cases}
\{1\}, & t>0,\\
[-1,1], & t=0,\\
\{-1\}, & t<0.
\end{cases}
\]
For $\mf u=(u_1,\dots,u_K)\in H$, let
\[
\mathcal S(\mf u)
:=
\left\{
\bs s=(s_1,\dots,s_K):
s_i(x)\in\operatorname{Sign}(u_i(x))
\text{ for a.e. }x\in\Omega,\ i=1,\dots,K
\right\}.
\]
By the standard calculus rules for the generalized gradient, one has
\begin{equation}\label{eq:generalized-gradient-E}
\partial E(\mf u)
\subseteq
\left\{
\mf g_{\bs s}:\bs s\in\mathcal S(\mf u)
\right\},
\end{equation}
where $\mf g_{\bs s}\in H^*$ is given by
\[
\mf g_{\bs s}[\bs\varphi]
=
\sum_{i=1}^K
\int_\Omega
\left(
\nabla u_i\cdot\nabla\varphi_i
+\lambda_i u_i\varphi_i
-\mu_i|u_i|^{2^*-2}u_i\varphi_i
\right)\,dx
-\beta
\sum_{i=1}^K
\int_\Omega
s_i
\prod_{\substack{j=1\\j\neq i}}^K|u_j|
\,\varphi_i\,dx,
\]
for every $\bs\varphi=(\varphi_1,\dots,\varphi_K)\in H$.

Following \cite{Chang_1981}, we say that $\mf u\in H$ is a critical point
of $E$ if
\[
0\in\partial E(\mf u).
\]

The next observation shows that this notion of critical point gives
precisely a weak solution of the original system \eqref{eq:system 1}.

\begin{lemma}\label{lem:clarke critical point are solutions}
Assume $q=1$ and let $\mf u\in H$. If
\[
\mf g_{\bs s}=0 \quad \text{for some } \bs s\in \cS(\mf u),
\]
then $\mf u$ is a weak solution of \eqref{eq:system 1}. In particular, by \eqref{eq:generalized-gradient-E}, this is true if $\mf u$ is a critical point.
\end{lemma}

\begin{proof}
The argument follows the same lines as
\cite{Clapp_Saldaña_Schiera_2026}.
By \eqref{eq:generalized-gradient-E}, there exists
$\bs s=(s_1,\dots,s_K)\in\mathcal S(\mf u)$ such that, for every
$i=1,\dots,K$,
\begin{equation}\label{eq:nonsmooth-equation}
-\Delta u_i+\lambda_i u_i
=
\mu_i|u_i|^{2^*-2}u_i
+
\beta s_i\prod_{j\neq i}|u_j|
\qquad\text{in }\Omega
\end{equation}
in the weak sense. Since $K=2^*$ and $u_j\in L^{2^*}(\Omega)$, H\"older's inequality gives
\[
\prod_{j\neq i}|u_j|
\in
L^{\frac{2^*}{2^*-1}}(\Omega).
\]
The right-hand side of \eqref{eq:nonsmooth-equation} therefore belongs
to $L^{\frac{2^*}{2^*-1}}(\Omega)$, and elliptic regularity yields $u_i\in
W_{\mathrm{loc}}^{2,\frac{2^*}{2^*-1}}(\Omega)$.
Therefore, we can apply \cite[Chapter~9, Remark~4, p.~314]{Brezis_2011}
to obtain
\[
\nabla u_i=0
\qquad\text{a.e. on }\{u_i=0\}.
\]
Applying the same result to each weak derivative $\partial_j u_i$,
$j=1,\dots,N$, we further obtain
\[
D^2u_i=0
\qquad\text{a.e. on }\{\nabla u_i=0\},
\]
and therefore
\[
\Delta u_i=0
\qquad\text{a.e. on }\{u_i=0\}.
\]
Hence, restricting \eqref{eq:nonsmooth-equation} to $\{u_i=0\}$,
we obtain
\[
s_i\prod_{j\neq i}|u_j|=0
\qquad\text{a.e. on }\{u_i=0\},
\]
since $\beta\neq0$. On the other hand,
$s_i=\operatorname{sign}(u_i)$ a.e. on $\{u_i\neq0\}$ by definition.
Consequently,
\[
s_i\prod_{j\neq i}|u_j|
=
\operatorname{sign}(u_i)
\prod_{j\neq i}|u_j|
\qquad\text{a.e. in }\Omega.
\]
Thus \eqref{eq:nonsmooth-equation} becomes
\[
-\Delta u_i+\lambda_i u_i
=
\mu_i|u_i|^{2^*-2}u_i
+
\beta\operatorname{sign}(u_i)
\prod_{j\neq i}|u_j|,
\qquad i=1,\dots,K,
\]
which is precisely \eqref{eq:system 1}.
\end{proof}

\section{Least energy solutions on bounded domains}\label{sec:least energy competitive}

The aim of this section is to prove Theorems \ref{thm:least energy competitive N=5,6} and \ref{thm:least energy competitive N=3,4}. Since the proof strategy is analogous, we present them in a unified way. 

Recalling the definition of $\cB$ in \eqref{eq:B definition}, since
\[
\int_\Omega\left(|\nabla u|^2+\lambda_i u^2\right)dx
\geq
\left(1+\frac{\lambda_i}{\lambda_1(\Omega)}\right)
\int_\Omega|\nabla u|^2dx,
\qquad i=1,\dots,K,
\]
it is standard to see that $\cB>0$. As pointed out in the introduction, by
\cite{BN} the scalar Brezis-Nirenberg problem \eqref{Brezis-Nirenberg} has a positive least energy solution
$u_{\mu_i}\in C^2(\Omega)\cap C(\overline{\Omega})$ whenever \eqref{eq:lambda-general} holds, with energy
\begin{equation}\label{eq:B_mu less than S}
	B_{\mu_i}:=
	\frac12\int_\Omega
	\left(
	|\nabla u_{\mu_i}|^2
	+\lambda_i u_{\mu_i}^2
	\right)dx
	-
	\frac{\mu_i}{2^*}
	\int_\Omega |u_{\mu_i}|^{2^*}dx
	<\frac1N
	\mu_i^{\frac{2-N}{2}}S^{N/2},
	\quad i=1,\dots,K.
\end{equation}

The strategy is to minimize the energy on $\cM$ while ruling out concentration phenomena. As anticipated above, we shall prove in Proposition \ref{prop: min sotto E} that concentration is excluded whenever the energy estimate \eqref{eq:upper bound B} holds:
\[
	\cB<\min\left\{
		\min_{1\leq \ell\leq K}\Big(\sum_{k\not=\ell} B_{\mu_k}
		+\frac1N\mu_\ell^{-\frac{N-2}{2}}S^{\frac N2}\Big),\cA
		\right\}.
\]
In the next subsection, we discuss how to obtain it in low and high dimensions; afterwards, we proceed with the proof of the main theorems.

\subsection{Energy estimate \eqref{eq:upper bound B}}

We first discuss the estimate in dimension $N=5,6$.
	
	\begin{lemma}\label{lem:upper-bound-B-N56}
		Let $K=3$, $\beta<0$, and $N=5,6$. Then \eqref{eq:upper bound B} holds.		
	\end{lemma}
	
	\begin{proof}
		We prove that
		\begin{equation}\label{eq:B-upper-12-3}
			\cB<
			B_{\mu_1}+B_{\mu_2}
			+\frac1N\mu_3^{-\frac{N-2}{2}}S^{\frac N2},
		\end{equation}
		the other two estimates following by permuting the indices. The last estimate in \eqref{eq:upper bound B} follows from \eqref{eq:A value beta negative} together with \eqref{eq:B_mu less than S}.
		
		Let $z_1,z_2$ be least energy positive solutions of
		\[
		\begin{cases}
			-\Delta z_i+\lambda_i z_i
			=\mu_i z_i^{2^*-1}
			&\text{in }\Omega,\\
			z_i=0
			&\text{on }\partial\Omega,
		\end{cases}
		\qquad i=1,2.
		\]
		Since $z_1,z_2\in C^1(\overline\Omega)$, for every sufficiently
		small $r>0$ we can choose $y_r\in\Omega$ such that
		$B_{2r}(y_r)\Subset\Omega$ and
		\begin{equation}\label{eq:zi-small-ball}
			z_1+z_2\leq Cr
			\qquad\text{in }B_{2r}(y_r),
		\end{equation}
		for some $C>0$ independent of $r$. Let now $\eta_r\in C_c^\infty(B_{2r}(y_r))$ be such that
		\[
		0\leq\eta_r\leq1,
		\qquad
		\eta_r\equiv1\quad\text{in }B_r(y_r),
		\qquad
		|\nabla\eta_r|\leq\frac Cr.
		\]
		Consider the Aubin--Talenti bubble $U_{\varepsilon,y_r}$, 
		and set	$w_{\varepsilon,r}:=\eta_rU_{\varepsilon,y_r}$.
		Recall that
		\[
		\int_{\mathbb R^N}U_{\varepsilon,y_r}^{2^*}\,dx
		=
		\int_{\mathbb R^N}|\nabla U_{\varepsilon,y_r}|^2\,dx
		=
		S^{N/2}.
		\]
		
		We first derive some estimates for the truncated bubble.
		For the $L^{2^*}$-norm, since
		$\eta_r\equiv1$ in $B_r(y_r)$ and $0\leq\eta_r\leq1$,
		\[
		\begin{aligned}
			0&\leq S^{N/2}-\int_\Omega w_{\varepsilon,r}^{2^*}\,dx
			=\int_{\mathbb R^N}(1-\eta_r^{2^*})U_{\varepsilon,y_r}^{2^*}\,dx
			\leq\int_{\mathbb R^N\setminus B_r(y_r)}U_{\varepsilon,y_r}^{2^*}\,dx
			\\
			&=C_N\varepsilon^N\int_r^\infty
			\frac{\rho^{N-1}}{(\varepsilon^2+\rho^2)^N}\,d\rho
			\leq C_N\varepsilon^N\int_r^\infty\rho^{-N-1}\,d\rho
			=C_N\left(\frac{\varepsilon}{r}\right)^N.
		\end{aligned}
		\]
		Hence
		\begin{equation*}
			\int_\Omega w_{\varepsilon,r}^{2^*}\,dx
			=
			S^{N/2}
			+
			O\left(
			\left(\frac{\varepsilon}{r}\right)^N
			\right).
		\end{equation*}
		
		We next estimate the Dirichlet energy. Since
		\[
		\nabla w_{\varepsilon,r}
		=
		\eta_r\nabla U_{\varepsilon,y_r}
		+
		U_{\varepsilon,y_r}\nabla\eta_r,
		\]
		we have
		\begin{align*}
			\left|
			\int_\Omega|\nabla w_{\varepsilon,r}|^2\,dx
			-
			S^{N/2}
			\right|
			&\leq
			\int_{\mathbb R^N\setminus B_r(y_r)}
			|\nabla U_{\varepsilon,y_r}|^2\,dx\\
			&\quad
			+
			C\int_{B_{2r}(y_r)\setminus B_r(y_r)}
			U_{\varepsilon,y_r}
			|\nabla U_{\varepsilon,y_r}|
			|\nabla\eta_r|\,dx
			+
			C\int_{B_{2r}(y_r)\setminus B_r(y_r)}
			U_{\varepsilon,y_r}^2
			|\nabla\eta_r|^2\,dx.
		\end{align*}
		Moreover,
		\[
		|\nabla U_{\varepsilon,y_r}(x)|
		\leq
		C\varepsilon^{\frac{N-2}{2}}
		\frac{|x-y_r|}
		{(\varepsilon^2+|x-y_r|^2)^{N/2}}.
		\]
		Therefore,
		\begin{align*}
			\int_{\mathbb R^N\setminus B_r(y_r)}
			|\nabla U_{\varepsilon,y_r}|^2\,dx
			&\leq
			C\varepsilon^{N-2}
			\int_r^\infty
			\frac{\rho^{N+1}}
			{(\varepsilon^2+\rho^2)^N}\,d\rho\\
			&\leq
			C\varepsilon^{N-2}
			\int_r^\infty\rho^{-N+1}\,d\rho=
			C\left(\frac{\varepsilon}{r}\right)^{N-2}.
		\end{align*}
		Furthermore,
		\begin{equation*}
			\int_{B_{2r}(y_r)\setminus B_r(y_r)}
			U_{\varepsilon,y_r}
			|\nabla U_{\varepsilon,y_r}|
			|\nabla\eta_r|\,dx
			\leq
			\frac Cr\varepsilon^{N-2}
			\int_r^{2r}\rho^{-N+2}\,d\rho\leq
			C\left(\frac{\varepsilon}{r}\right)^{N-2},
		\end{equation*}
		and
		\begin{equation*}
			\int_{B_{2r}(y_r)\setminus B_r(y_r)}
			U_{\varepsilon,y_r}^2
			|\nabla\eta_r|^2\,dx
			\leq
			\frac C{r^2}\varepsilon^{N-2}
			\int_r^{2r}\rho^{-N+3}\,d\rho\\
			\leq
			C\left(\frac{\varepsilon}{r}\right)^{N-2}.
		\end{equation*}
		Thus
		\begin{equation*}
			\int_\Omega|\nabla w_{\varepsilon,r}|^2\,dx
			=
			S^{N/2}
			+
			O\left(
			\left(\frac{\varepsilon}{r}\right)^{N-2}
			\right).
		\end{equation*}
		
		We now estimate the $L^2$-norm. Since $N\geq5$,
		$U_{1,0}\in L^2(\mathbb R^N)$ and, by the change of variables
		$x-y_r=\varepsilon z$,
		\begin{align*}
			\int_{\mathbb R^N}
			U_{\varepsilon,y_r}^2\,dx
			&=
			\alpha_N^2\varepsilon^{N-2}
			\int_{\mathbb R^N}
			\frac{dx}
			{(\varepsilon^2+|x-y_r|^2)^{N-2}}\\
			&=
			\varepsilon^2
			\alpha_N^2
			\int_{\mathbb R^N}
			\frac{dz}{(1+|z|^2)^{N-2}}=
			C_N\varepsilon^2,
		\end{align*}
		where
		\[
		C_N:=\int_{\mathbb R^N}U_{1,0}^2\,dx>0.
		\]
		Therefore,
		\begin{align*}
			0
			&\leq
			C_N\varepsilon^2
			-
			\int_\Omega w_{\varepsilon,r}^2\,dx
			\leq
			\int_{\mathbb R^N\setminus B_r(y_r)}
			U_{\varepsilon,y_r}^2\,dx\\
			&\leq
			C\varepsilon^{N-2}
			\int_r^\infty
			\frac{\rho^{N-1}}
			{(\varepsilon^2+\rho^2)^{N-2}}\,d\rho
			\leq
			C\varepsilon^{N-2}
			\int_r^\infty\rho^{-N+3}\,d\rho=
			C\frac{\varepsilon^{N-2}}{r^{N-4}}.
		\end{align*}
		Hence
		\begin{equation*}
			\int_\Omega w_{\varepsilon,r}^2\,dx
			=
			C_N\varepsilon^2
			+
			O\left(
			\frac{\varepsilon^{N-2}}{r^{N-4}}
			\right).
		\end{equation*}

		Finally, recalling that
		\[
		q=\frac{2^*}{3}
		=
		\frac{2N}{3(N-2)},
		\qquad
		\frac{N-2}{2}q=\frac N3,
		\]
		we obtain
		\begin{align*}
			\int_{B_{2r}(y_r)}
			w_{\varepsilon,r}^q\,dx
			&\leq
			C\varepsilon^{N/3}
			\int_0^{2r}
			\frac{\rho^{N-1}}
			{(\varepsilon^2+\rho^2)^{N/3}}\,d\rho\\
			&\leq
			C\varepsilon^{N/3}
			\int_0^{2r}\rho^{\frac N3-1}\,d\rho \leq
			C\varepsilon^{N/3}r^{N/3}.
		\end{align*}
		Therefore
		\begin{equation}\label{eq:bubble-Lq}
			\int_{B_{2r}(y_r)}
			w_{\varepsilon,r}^q\,dx
			=
			O\left(
			\varepsilon^{N/3}r^{N/3}
			\right).
		\end{equation}
		
		From now on, we choose
		\[
		r=r(\varepsilon):=\varepsilon^{1/4}.
		\]
		Then
		\[
		\left(\frac{\varepsilon}{r}\right)^{N-2}
		=
		\begin{cases}
			\varepsilon^{9/4},&N=5,\\
			\varepsilon^3,&N=6,
		\end{cases}
		\]
		and
		\[
		\frac{\varepsilon^{N-2}}{r^{N-4}}
		=
		\begin{cases}
			\varepsilon^{11/4},&N=5,\\
			\varepsilon^{7/2},&N=6,
		\end{cases}
		\]
		while
		\[
		\left(\frac{\varepsilon}{r}\right)^N
		=
		\varepsilon^{3N/4}.
		\]
		In particular, all these quantities are $o(\varepsilon^2)$.
		Writing $w_\varepsilon:=w_{\varepsilon,r(\varepsilon)}$,
		we therefore have
		\begin{equation}\label{eq:bubble-estimates-final}
			\begin{aligned}
				\int_\Omega w_\varepsilon^{2^*}\,dx
				&=
				S^{N/2}+o(\varepsilon^2),\\
				\int_\Omega|\nabla w_\varepsilon|^2\,dx
				&=
				S^{N/2}+o(\varepsilon^2),\\
				\int_\Omega w_\varepsilon^2\,dx
				&=
				C_N\varepsilon^2+o(\varepsilon^2),
			\end{aligned}
		\end{equation}
		and hence
		\begin{equation*}
			\int_\Omega	\left(|\nabla w_{\varepsilon}|^2 +\lambda_3w_{\varepsilon}^2	\right)\,dx	=		S^{N/2}	+\lambda_3C_N\varepsilon^2	+		o(\varepsilon^2).
		\end{equation*}
		
		Consider now the scalar fiber map
		\[
		I_3(s)
		:=
		\frac{s^2}{2}
		\int_\Omega	\left(|\nabla w_{\varepsilon}|^2 +\lambda_3w_{\varepsilon}^2\right)\,dx
		-
		\frac{\mu_3s^{2^*}}{2^*}
		\|w_\varepsilon\|_{2^*}^{2^*}.
		\]
		Its maximum is attained at
		\begin{equation}\label{eq:s-epsilon}
			s_\varepsilon
			=
			\left(
			\frac{
				\int_\Omega	\left(|\nabla w_{\varepsilon}|^2 +\lambda_3w_{\varepsilon}^2\right)\,dx
			}{
				\mu_3\|w_\varepsilon\|_{2^*}^{2^*}
			}
			\right)^{\frac1{2^*-2}},
		\end{equation}
		and therefore
		\begin{align}
			\max_{s>0}I_3(s)
			&=
			\frac1N
			\frac{
				\left(
				\int_\Omega	\left(|\nabla w_{\varepsilon}|^2 +\lambda_3w_{\varepsilon}^2\right)\,dx
				\right)^{\frac{2^*}{2^*-2}}
			}{
				\left(
				\mu_3
				\|w_\varepsilon\|_{2^*}^{2^*}
				\right)^{\frac2{2^*-2}}
			}
			\notag\\
			&=
			\frac1N
			\mu_3^{-\frac{N-2}{2}}
			S^{N/2}
			\left(
			1+\frac N2\lambda_3C_NS^{-N/2}\varepsilon^2
			+o(\varepsilon^2)
			\right).
			\label{eq:max-I3-expansion}
		\end{align}
		Since $\lambda_3<0$, it follows that
		\begin{equation*}
			\max_{s>0}I_3(s)
			<
			\frac1N
			\mu_3^{-\frac{N-2}{2}}
			S^{N/2}
		\end{equation*}
		for $\varepsilon>0$ sufficiently small.
		Moreover, by \eqref{eq:s-epsilon} and
		\eqref{eq:bubble-estimates-final},
		\[
		s_\varepsilon
		\longrightarrow
		\mu_3^{-\frac{N-2}{4}}
		\qquad\text{as }\varepsilon\to0^+,
		\]
		and hence $s_\varepsilon$ is bounded from above and away from zero.
		
		We next estimate the interaction term. By
		\eqref{eq:zi-small-ball}, \eqref{eq:bubble-Lq}, and the fact that
		$\supp w_\varepsilon\subset B_{2r}(y_r)$,
		\begin{equation*}
			\int_\Omega
			z_1^qz_2^qw_\varepsilon^q\,dx
			\leq
			Cr^{2q}
			\int_{B_{2r}(y_r)}w_\varepsilon^q\,dx\leq
			C\varepsilon^{N/3}r^{N/3+2q}.
		\end{equation*}
		With $r=\varepsilon^{1/4}$, this becomes
		\[
		\int_\Omega
		z_1^qz_2^qw_\varepsilon^q\,dx
		=
		\begin{cases}
			O(\varepsilon^{5/2}),&N=5,\\
			O(\varepsilon^{3}),&N=6,
		\end{cases}
		\]
		and hence
		\begin{equation}\label{eq:interaction-small}
			\int_\Omega
			z_1^qz_2^qw_\varepsilon^q\,dx
			=
			o(\varepsilon^2).
		\end{equation}
		
		We finally project the competitor onto $\mathcal M$.
		We claim that, for $\varepsilon>0$ sufficiently small, there exist
		$t_1^\varepsilon,t_2^\varepsilon,t_3^\varepsilon>0$ such that
		\begin{equation}\label{eq:competitor-M}
			(t_1^\varepsilon z_1,
			t_2^\varepsilon z_2,
			t_3^\varepsilon s_\varepsilon w_\varepsilon)
			\in\mathcal M.
		\end{equation}
		Set
		\[
		K_\varepsilon
		:=
		\int_\Omega
		z_1^qz_2^q(s_\varepsilon w_\varepsilon)^q\,dx.
		\]
		By \eqref{eq:interaction-small} and the boundedness of
		$s_\varepsilon$, $K_\varepsilon=o(\varepsilon^2)$.
		
Since $z_1,z_2$ and $s_\varepsilon w_\varepsilon$ belong to their
respective scalar Nehari manifolds, we have
\[
B_{\mu_i}=\frac1N\|z_i\|_{\lambda_i}^2,
\qquad i=1,2,
\]
and
\[
\frac1N\|s_\varepsilon w_\varepsilon\|_{\lambda_3}^2
=
\max_{s>0}I_3(s).
\]

Condition \eqref{eq:competitor-M} is equivalent to
\begin{equation*}
\begin{cases}
\bigl((t_i^\varepsilon)^2-(t_i^\varepsilon)^{2^*}\bigr)
\|z_i\|_{\lambda_i}^2
-\beta(t_1^\varepsilon t_2^\varepsilon t_3^\varepsilon)^q
K_\varepsilon=0,
& i=1,2,\\[2mm]
\bigl((t_3^\varepsilon)^2-(t_3^\varepsilon)^{2^*}\bigr)
\|s_\varepsilon w_\varepsilon\|_{\lambda_3}^2
-\beta(t_1^\varepsilon t_2^\varepsilon t_3^\varepsilon)^q
K_\varepsilon=0.
\end{cases}
\end{equation*}
Define
\[
\Phi_i(\bs t,K)
:=
(t_i^2-t_i^{2^*})\|z_i\|_{\lambda_i}^2
-\beta(t_1t_2t_3)^qK,
\qquad i=1,2,
\]
and
\[
\Phi_3(\bs t,K)
:=
(t_3^2-t_3^{2^*})
\|s_\varepsilon w_\varepsilon\|_{\lambda_3}^2
-\beta(t_1t_2t_3)^qK.
\]
Then $\Phi(1,1,1,0)=0$ and
\[
D_{\bs t}\Phi(1,1,1,0)
=
(2-2^*)
\begin{pmatrix}
\|z_1\|_{\lambda_1}^2&0&0\\
0&\|z_2\|_{\lambda_2}^2&0\\
0&0&\|s_\varepsilon w_\varepsilon\|_{\lambda_3}^2
\end{pmatrix}.
\]
Since $\|z_1\|_{\lambda_1}^2$ and
$\|z_2\|_{\lambda_2}^2$ are fixed and positive, while
$\|s_\varepsilon w_\varepsilon\|_{\lambda_3}^2$
is bounded from above and away from zero for
$\varepsilon>0$ sufficiently small, this matrix
is uniformly invertible. Therefore, by the
implicit function theorem, there exist
$\delta>0$ and $\varepsilon_0>0$ such that,
for every $\varepsilon\in(0,\varepsilon_0)$,
there exists a unique $C^1$ function
$\bs\ell_\varepsilon:(-\delta,\delta)\to\mathbb R^3$
such that $\bs\ell_\varepsilon(0)=(1,1,1)$ and
$\Phi(\bs t,K)=0$ with
$(\bs t,K)\in B_\delta(1,1,1)\times(-\delta,\delta)$
if and only if $\bs t=\bs\ell_\varepsilon(K)$.

Since $K_\varepsilon=o(\varepsilon^2)$,
we have $K_\varepsilon\in(-\delta,\delta)$
for $\varepsilon>0$ sufficiently small.
Hence there exists a unique
\[
\bs t^\varepsilon
=
(t_1^\varepsilon,t_2^\varepsilon,t_3^\varepsilon)
=
\bs\ell_\varepsilon(K_\varepsilon)
\in B_\delta(1,1,1)
\]
such that $\Phi(\bs t^\varepsilon,K_\varepsilon)=0$.
Equivalently, $(t_1^\varepsilon z_1,
t_2^\varepsilon z_2,
t_3^\varepsilon s_\varepsilon w_\varepsilon)
\in\mathcal M$.

Moreover, differentiating the implicit equation
with respect to $K$ at $K=0$, we obtain
\[
\bs\ell_\varepsilon'(0)
=
-\frac{\beta}{2^*-2}
\begin{pmatrix}
\|z_1\|_{\lambda_1}^{-2}\\
\|z_2\|_{\lambda_2}^{-2}\\
\|s_\varepsilon w_\varepsilon\|_{\lambda_3}^{-2}
\end{pmatrix}.
\]
Consequently,
\begin{equation*}
t_i^\varepsilon
=
1+O(|\beta|K_\varepsilon)
=
1+o(\varepsilon^2),
\qquad i=1,2,3.
\end{equation*}
We conclude that
\begin{align*}
\cB
&\leq
E(t_1^\varepsilon z_1,
t_2^\varepsilon z_2,
t_3^\varepsilon s_\varepsilon w_\varepsilon)\\
&=
\frac1N\left(
\sum_{i=1}^2
(t_i^\varepsilon)^2\|z_i\|_{\lambda_i}^2
+
(t_3^\varepsilon)^2
\|s_\varepsilon w_\varepsilon\|_{\lambda_3}^2
\right)\\
&=
B_{\mu_1}+B_{\mu_2}
+\frac1N\|s_\varepsilon w_\varepsilon\|_{\lambda_3}^2
+o(\varepsilon^2).
\end{align*}
By \eqref{eq:max-I3-expansion},
\[
\frac1N\|s_\varepsilon w_\varepsilon\|_{\lambda_3}^2
=
\frac1N\mu_3^{-\frac{N-2}{2}}S^{N/2}
+
\frac12\lambda_3 C_N
\mu_3^{-\frac{N-2}{2}}\varepsilon^2
+
o(\varepsilon^2).
\]
Therefore,
\[
\cB
\leq
B_{\mu_1}+B_{\mu_2}
+\frac1N\mu_3^{-\frac{N-2}{2}}S^{N/2}
-C\varepsilon^2
+o(\varepsilon^2)
\]
for some $C>0$, since $\lambda_3<0$.
Hence, for $\varepsilon>0$ sufficiently small,
the desired estimate \eqref{eq:B-upper-12-3} follows.
	\end{proof}
	
	\begin{remark}
		Notice that the previous proof also applies when $N=6$, in which
		case $q=1$. Indeed, only the explicit component-wise Nehari
		identities are used, and no differentiability of the interaction
		functional on $H$ is required.
	\end{remark}

When $N=3,4$ and $\beta<0$, the Sobolev exponent is $2^* = 6$ (for $N=3$) or $2^* = 4$ (for $N=4$), thus the competitor used in dimension $N=5,6$ cannot be used because the Talenti bubble is not in $L^2$. However, using the second Dirichlet eigenfunction $e_2$ of $-\Delta$ on $\Omega$, we can establish an abstract condition on $\Omega$ that guarantees the energy upper bound \eqref{eq:upper bound B}.

For any $\sigma \in H_0^1(\Omega) \setminus \{0\}$ and $i \in \{1,\dots,K\}$, we define the Rayleigh-type functional
\[
R_i(\sigma) := \max_{t>0} I_i(t\sigma) = \frac{1}{N} \mu_i^{-\frac{N-2}{2}} \frac{\|\sigma\|_{\lambda_i}^N}{\|\sigma\|_{2^*}^{N}},
\]
where $\|\sigma\|_{\lambda_i}^2 := \int_\Omega (|\nabla \sigma|^2 + \lambda_i \sigma^2) \, dx$ and $I_i(u) := \frac{1}{2}\|u\|_{\lambda_i}^2 - \frac{\mu_i}{2^*}\|u\|_{2^*}^{2^*}$.
The first lemma is a general observation.

\begin{lemma}\label{lemma:general B bound with raileigh}
Let  $\sigma_i, \sigma_j \in H_0^1(\Omega) \setminus \{0\}$ be two functions such that
\[
\sigma_i \sigma_j \equiv 0 \qquad \text{a.e. in $\Omega$},
\]
where $i,j\in \{1,\dots,K\}$. Then
\[
\cB \le R_i(\sigma_i) + R_j(\sigma_j) + \sum_{k\not=i,j}B_{\mu_k}.
\]
\end{lemma}

\begin{proof}
Let $z_k \in H_0^1(\Omega)$ with $k\in\{1,\dots,K\}\backslash\{i,j\}$ be a positive least energy solution of the scalar problem
\begin{equation}\label{scalar problem}
-\Delta z_k + \lambda_k z_k = \mu_k z_k^{2^*-1} \quad \text{in } \Omega, \qquad z_k \in H_0^1(\Omega).
\end{equation}
Since $\sigma_i \sigma_j =0$ a.e., the energy functional evaluated on $(t_i \sigma_i, t_j \sigma_j, (z_k)_{k\not=i,j})$ completely decouples:
\[
E(t_i \sigma_i, t_j \sigma_j, (z_k)_{k\not=i,j}) = I_i(t_i \sigma_i) + I_j(t_j \sigma_j) + \sum_{k\not=i,j}I_k(z_k).
\]
Moreover, choosing 
$$
t_i = t_{\sigma_i} := \left( \frac{\|\sigma_i\|_{\lambda_i}^2}{\mu_i\|\sigma_i\|_{2^*}^{2^*}} \right)^{\frac{1}{2^*-2}}
$$
and $t_j$ analogously, recalling that $z_k$ already belongs to its scalar Nehari manifold, we obtain
\[
(t_{\sigma_i} \sigma_i, t_{\sigma_j} \sigma_j, (z_k)_{k\not=i,j}) \in \cM.
\]
By definition of $\cB = \inf_{\cM} E$, it follows that
\[
\cB \le E(t_{\sigma_i} \sigma_i, t_{\sigma_j} \sigma_j, (z_k)_{k\not=i,j}) = I_i(t_{\sigma_i} \sigma_i) + I_j(t_{\sigma_j} \sigma_j) + \sum_{k\not=i,j}I_k(z_k) = R_i(\sigma_i) + R_j(\sigma_j) + \sum_{k\not=i,j}B_{\mu_k}.
\]
\end{proof}

Let $e_2 \in H_0^1(\Omega) \cap C^1(\overline{\Omega})$ be a second eigenfunction of $-\Delta$ on $\Omega$, corresponding to $\Lambda_2 = \Lambda_2(\Omega)$, i.e.,
\[
\begin{cases}
    -\Delta e_2 = \Lambda_2 e_2 & \text{in } \Omega,\\
    e_2 = 0 & \text{on } \partial\Omega.
\end{cases}
\]

Set $e_2^+ := \max\{e_2,0\}$ and
$e_2^- := \max\{-e_2,0\}$. Note that
$e_2^+,e_2^- \not\equiv 0$ and have disjoint positivity sets
\[
\Omega^+ := \{x\in\Omega : e_2(x)>0\},
\qquad
\Omega^- := \{x\in\Omega : e_2(x)<0\}.
\]
We also introduce the notation
\[
\overline{\mu}
:= \max_{1\leq i\leq K}\mu_i^{-\frac{N-2}{2}},
\qquad
\underline{\mu}
:= \min_{1\leq i\leq K}\mu_i^{-\frac{N-2}{2}},
\qquad
\overline{\lambda}
:= \max_{1\leq i\leq K}\lambda_i.
\]
We can now establish the following stronger estimate,
which is uniform with respect to $\beta$.
To emphasize the dependence of the energy level
$\cB$ on $\beta$, we denote it by $\cB_\beta$.

\begin{lemma}\label{lem:abstract-bound-N34}
Let $K\geq 3$ and $N\in\{3,4\}$.
Suppose that $(\Lambda)$ holds and that
$\Omega$ satisfies condition~\eqref{eq:condition-C}.
Then
\[
\sup_{\beta<0}\cB_\beta
<
\min\left\{
\min_{1\leq \ell\leq K}
\left(
\sum_{\substack{k=1\\k\neq\ell}}^K B_{\mu_k}
+\frac1N\mu_\ell^{-\frac{N-2}{2}}S^{N/2}
\right),
\cA
\right\}.
\]
In particular, \eqref{eq:upper bound B}
holds for every $\beta<0$.
\end{lemma}

\begin{proof}
We may assume, as in the statement of Theorem \ref{thm:least energy competitive N=3,4}, that $B_{\mu_1}\leq \dots \leq B_{\mu_K}$.
For any two distinct indices $i,j\in  \{1,\dots,K\}$. Applying Lemma \ref{lemma:general B bound with raileigh} with
\[
\sigma_i=e_2^+,~\sigma_j=e_2^-,
\]
and with $z_k$ the positive least energy solution of the scalar problem  \eqref{scalar problem} for $k\not=i,j$, we obtain 
\begin{equation}\label{Bijk}
    \cB_\beta\leq R_i(\sigma_i) + R_j(\sigma_j) + \sum_{k\not=i,j}B_{\mu_k}.
\end{equation}
Testing $-\Delta e_2 = \Lambda_2 e_2$ with $e_2^+$ and $e_2^-$ respectively, we have
\[
\int_\Omega |\nabla e_2^\pm|^2 \, dx + \lambda_i \int_\Omega (e_2^\pm)^2 \, dx = (\Lambda_2 + \lambda_i) \int_\Omega (e_2^\pm)^2 \, dx.
\]
By Hölder's inequality on $\Omega^\pm$,
\[
\int_{\Omega^\pm} (e_2^\pm)^2 \, dx \leq \left( \int_{\Omega^\pm} (e_2^\pm)^{2^*} \, dx \right)^{\frac{2}{2^*}} |\Omega^\pm|^{\frac{2}{N}},
\]
and noting that $2/2^* = (N-2)/N$, we obtain
\[
R_i(e_2^+) \leq \frac{1}{N} \mu_i^{-\frac{N-2}{2}} (\Lambda_2 + \lambda_i)^{N/2} |\Omega^+|,
\]
and similarly
\[
R_j(e_2^-) \leq \frac{1}{N} \mu_j^{-\frac{N-2}{2}} (\Lambda_2 + \lambda_j)^{N/2} |\Omega^-|.
\]
Summing these two inequalities and using $|\Omega^+| + |\Omega^-| \leq |\Omega|$, we find
\[
R_i(e_2^+) + R_j(e_2^-) \leq \frac{1}{N} \overline{\mu} \left( \Lambda_2(\Omega) + \overline{\lambda} \right)^{N/2} |\Omega|,
\]
which together with  \eqref{Bijk} shows that
\[
\cB_\beta\leq \frac{1}{N} \overline{\mu} \left( \Lambda_2(\Omega) + \overline{\lambda} \right)^{N/2} |\Omega|+ \sum_{k\not=i,j}B_{\mu_k}.
\]
To minimize the right hand side, we choose $i=K-1, j=K$. Therefore,
\[
\begin{aligned}
\sup_{\beta<0}B_\beta
&\le
\frac{\overline\mu}{N}
\bigl(\Lambda_2(\Omega)+\overline\lambda\bigr)^{N/2}
|\Omega|
+\sum_{i=1}^{K-2} B_{\mu_i}<
\sum_{i=1}^{K-1} B_{\mu_i}
+\frac{\underline\mu}{N}S^{N/2}\\
&\le
\min_{1\leq\ell\leq K}
\left(
\sum_{i\ne\ell}B_{\mu_i}
+\frac1N\mu_\ell^{-\frac{N-2}{2}}S^{N/2}
\right),
\end{aligned}
\]
where the strict inequality follows from condition \eqref{eq:condition-C}. The last estimate follows again from \eqref{eq:A value beta negative} together with \eqref{eq:B_mu less than S}, and this concludes the proof.
\end{proof}

\subsection{General tools for the minimization on $\cM$}

We collect below some preliminary results concerning the component-wise
Nehari manifold. They are direct adaptations to the critical case
$Kq=2^*$, with $\beta<0$, of the corresponding results proved in
\cite{Giaretto_Soave_2026}. Since the proofs are identical, we only
recall the main arguments when useful.

For $i=1,\dots,K$, define
\[
G_{\beta,i}(\mf u)
:=
\|u_i\|_{\lambda_i}^2
-\mu_i |u_i|_{2^*}^{2^*}
-\beta\int_\Omega\prod_{j=1}^K |u_j|^q\,dx,
\]
so that
\[
\mathcal M
=
\left\{
\mf u\in H:
u_i\not\equiv0,\ 
G_{\beta,i}(\mf u)=0
\quad\text{for every }i=1,\dots,K
\right\}.
\]
For $\mf u\in\mathcal M$ and $j=1,\dots,K$, let
\[
(\bs\varphi_j(\mf u))_k
=
\begin{cases}
u_j, & k=j,\\
0, & k\neq j,
\end{cases}
\qquad k=1,\dots,K.
\]
If $q>1$, we define $\bs M_\beta(\mf u)\in\R^{K\times K}$ by
\[
(\bs M_\beta(\mf u))_{ij}
:=
G_{\beta,i}'(\mf u)[\bs\varphi_j(\mf u)].
\]
If $q=1$, for each $i=1,\dots,K$ we choose an arbitrary
\[
\xi_i\in\partial G_{\beta,i}(\mf u)
\]
and define
\[
(\bs M_\beta(\mf u))_{ij}
:=
\langle \xi_i,\bs\varphi_j(\mf u)\rangle.
\]
This definition is independent of the choice of $\xi_i$. Indeed, the
contribution of the interaction term to any
$\xi_i\in\partial G_{\beta,i}(\mf u)$ has the form
\[
-\beta\int_\Omega
s_i\prod_{k\neq i}|u_k|\,\varphi_i\,dx
-\beta\sum_{\ell\neq i}
\int_\Omega
s_\ell |u_i|
\prod_{k\neq i,\ell}|u_k|\,\varphi_\ell\,dx,
\]
where
\[
s_k(x)\in\operatorname{Sign}(u_k(x))
\qquad\text{for a.e. }x\in\Omega.
\]
Evaluating at $\bs\varphi_j(\mf u)$ and observing that
\[
s_j u_j=|u_j|
\qquad\text{a.e. in }\Omega,
\]
regardless of the choice of $s_j$ on $\{u_j=0\}$, the contribution of
the interaction term is always
\[
-\beta\int_\Omega\prod_{k=1}^K|u_k|\,dx.
\]
Therefore, in both cases $q>1$ and $q=1$, using $Kq=2^*$ and
$G_{\beta,i}(\mf u)=0$, we have
\[
(\bs M_\beta(\mf u))_{ij}
=
\begin{cases}
(2-2^*)\mu_i|u_i|_{2^*}^{2^*}
+(2-q)\beta\displaystyle\int_\Omega
\prod_{k=1}^K|u_k|^q\,dx,
& i=j,\\[0.5em]
-q\beta\displaystyle\int_\Omega
\prod_{k=1}^K|u_k|^q\,dx,
& i\neq j.
\end{cases}
\]
For further use in the case $q=1$, for every $i=1,\dots,K$ and
$\bs t=(t_1,\dots,t_K)\in\cS(\mf u)$, we denote by
$G_{\beta,i,\bs t}'(\mf u)\in H'$ the functional defined by
\begin{align*}
G_{\beta,i,\bs t}'(\mf u)[\bs\varphi]
&:=
2\int_\Omega
\left(
\nabla u_i\cdot\nabla\varphi_i
+\lambda_i u_i\varphi_i
\right)\,dx
-
2^*\mu_i
\int_\Omega
|u_i|^{2^*-2}u_i\varphi_i\,dx
\\
&\quad
-
\beta
\sum_{j=1}^K
\int_\Omega
t_j
\prod_{\ell\neq j}|u_\ell|
\,\varphi_j\,dx,
\end{align*}
for every $\bs\varphi\in H$. 
With this notation, for every $i=1,\dots,K$,
\[
\partial G_{\beta,i}(\mf u)
\subseteq
\left\{
G_{\beta,i,\bs t}'(\mf u):
\bs t\in\cS(\mf u)
\right\}.
\]

\begin{lemma}\label{lem:uniform-bounds-components}
		Let $\beta<0$. Then there exist constants $C_1,C_2>0$ such that
		\[
		C_1\leq |u_i|_{2^*}\leq C_2
		\qquad
		\text{for every } \mf u\in\mathcal M
		\text{ such that } E(\mf u)\leq \cA,
		\]
		for every $i=1,\ldots,K$.
	\end{lemma}
	
	\begin{proof}
		The proof is identical to that of \cite[Lemma 4.5]{Giaretto_Soave_2026}.
	\end{proof}

\begin{lemma}\label{lemma:M_beta_negative}
Let $\beta<0$. Then
$\bs M_\beta(\mf u)$ is invertible, uniformly with respect to $\mf u\in \cM$ with $E(\mf u)\leq \cA$.
\end{lemma}

\begin{proof}
The proof is identical to that of
\cite[Lemma~4.3]{Giaretto_Soave_2026}.
Indeed, since $\mf u\in\mathcal M$ and $\beta<0$,
\[
-(M_\b(\mf u))_{ii}-\sum_{j\neq i}|(M_\b(\mf u))_{ij}|=(2^*-2)\norm{u_{i}}_{\lambda_i}^2\geq c>0,
\]
for some $c>0$ independent of $\mf u$, where the last inequality follows by Lemma \ref{lem:uniform-bounds-components}.
The conclusion now follows from Gershgorin's circle theorem.
\end{proof}

\begin{proposition}\label{prop:minimizers_are_solutions}
Let $\beta<0$ and let $\mf u\in\mathcal M$ be a minimizer
of $E$ on $\mathcal M$. Then $\mf u$ is a critical point of $E$
and hence a weak solution of \eqref{eq:system 1}.
\end{proposition}

\begin{proof}
When $q>1$, the proof is identical to that of
\cite[Proposition~4.7]{Giaretto_Soave_2026}.
Arguing as in \cite[Lemma~4.4]{Giaretto_Soave_2026},
one can see that $\cM$ is a $C^1$ manifold, hence
the classical Lagrange multiplier rule applies,
while \Cref{lemma:M_beta_negative} implies that
all the Lagrange multipliers vanish.

When $q=1$, we use Clarke's generalized Lagrange
multiplier rule in \cite[Theorem~1 and Corollary~1]{Clarke_1976}.
Since $E$ and $G_{\beta,i}$ are locally Lipschitz,
there exist $\lambda_0\geq0$ and
$\lambda_1,\ldots,\lambda_K\in\R$, not all zero,
such that
\[
0\in
\lambda_0\partial E(\mf u)
+\sum_{i=1}^K\lambda_i\partial G_{\beta,i}(\mf u).
\]
Thus, there exist $\xi\in\partial E(\mf u)$ and
$\xi_i\in\partial G_{\beta,i}(\mf u)$ such that
\[
\lambda_0\xi+\sum_{i=1}^K\lambda_i\xi_i=0.
\]
Testing against $\bs\varphi_j(\mf u)$,
$j=1,\ldots,K$, and recalling that
\[
\langle\xi,\bs\varphi_j(\mf u)\rangle
=G_{\beta,j}(\mf u)=0,
\]
we obtain
\[
\bs M_\beta(\mf u)^T\bs\lambda=\bs0,
\qquad
\bs\lambda=(\lambda_1,\ldots,\lambda_K).
\]
By \Cref{lemma:M_beta_negative},
$\bs M_\beta(\mf u)$ is invertible, hence
$\bs\lambda=\bs0$. Since the multipliers
are not all zero, necessarily $\lambda_0>0$,
and therefore
\[
0\in\partial E(\mf u).
\]
By Lemma \ref{lem:clarke critical point are solutions}, we conclude.
\end{proof}

We will make use of the following Brezis-Lieb type lemma.
\begin{lemma}\label{lem:BL-interaction}
	Assume $q\geq1$. Let $\mf u_n,\mf u\in H$ be such that
	\[
	\mf u_n\rightharpoonup \mf u
	\qquad\text{in }H.
	\]
	Then, up to subsequences,
	\[
	\int_\Omega
	\left(
	\prod_{i=1}^K |u_{i,n}|^q
	-
	\prod_{i=1}^K |u_{i,n}-u_i|^q
	-
	\prod_{i=1}^K |u_i|^q
	\right)\,dx
	\longrightarrow 0.
	\]
\end{lemma}

\begin{proof}
	See \cite[Lemma~4.4]{FaGiSo26}.
\end{proof}

We are now in a position to show that, under the energy estimate \eqref{eq:upper bound B}, the level $\cB$ is attained by a least energy fully non-trivial solution.

\begin{proposition}\label{prop: min sotto E}
    Let $K\geq 3$, $N\geq 3$ and $\beta<0$. If \eqref{eq:upper bound B} holds, $\cB$ is achieved by a least energy fully non-trivial solution of \eqref{eq:system 1}.
\end{proposition}
\begin{proof}
    \textit{Case $q>1$:} in this case, $E\in C^1(H)$. Recall that $E$ is coercive and bounded from below on $\mathcal M$. By Ekeland's variational principle, and using also that $\cM$ is a $C^1$ manifold as observed in Proposition \ref{prop:minimizers_are_solutions}, there exists
		$(\mf u_n)\subset\mathcal M$ such that
		\begin{equation*}
			E(\mf u_n)\to \cB,
			\qquad
			(E|_{\mathcal M})'(\mf u_n)\to0.
		\end{equation*}
		In particular, $(\mf u_n)$ is bounded in $H$ and,
		up to a subsequence, $\mf u_n\rightharpoonup\mf u$ weakly in $H$.
		Using Lemma \ref{lemma:M_beta_negative} and arguing as in \cite[Lemma 4.6]{Giaretto_Soave_2026}, one can
		show that $(\mf u_n)$ is actually an unconstrained Palais--Smale
		sequence, namely $E'(\mf u_n)\to0$ in $H'$. By weak convergence, the Sobolev embeddings, and standard
		weak-continuity arguments for the nonlinear terms, we obtain $E'(\mf u) = 0$.
        
Set
\[
	\mf w_n:=\mf u_n-\mf u.
\]
By the Brezis--Lieb lemma, for every $i=1,\ldots,K$,
\begin{equation*}
    \|u_{i,n}\|_{2^*}^{2^*}=\|u_i\|_{2^*}^{2^*}+\|w_{i,n}\|_{2^*}^{2^*}+o(1).
\end{equation*}
Moreover, Lemma~\ref{lem:BL-interaction} gives
		\[
		\int_\Omega
		\left(
		\prod_{i=1}^K |u_{i,n}|^q
		-
		\prod_{i=1}^K |w_{i,n}|^q
		-
		\prod_{i=1}^K |u_i|^q
		\right)\,dx
		=o(1).
		\]
		Since $w_{i,n}\to0$ strongly in $L^2(\Omega)$, combining these
		relations with $G_{i,\b}(\mf u_n)=0$, $E'(\mf u_n)[\mf u_n]\to0$ (resp. $\langle\xi_n,\mf u_n\rangle\to 0$) and $E'(\mf u)=0$ (resp. $\mf g_{\bs u}=0$) yields, for every
		$i=1,\ldots,K$,
		\begin{equation}\label{eq:remainder-nehari}
			\int_\Omega |\nabla w_{i,n}|^2dx
			-
			\mu_i\int_\Omega |w_{i,n}|^{2^*}dx
			-
			\beta\int_\Omega\prod_{\ell=1}^K|w_{\ell,n}|^qdx
			=o(1).
		\end{equation}
		Furthermore,
		\begin{equation}\label{eq:energy-splitting2}
			E(\mf u_n)
			=
			E(\mf u)+I(\mf w_n)+o(1),
		\end{equation}
		where
		\[
		I(\mf w)
		:=
		\frac12\sum_{i=1}^K\int_\Omega|\nabla w_i|^2dx
		-
		\frac1{2^*}\int_\Omega
		\left(
		\sum_{i=1}^K\mu_i|w_i|^{2^*}
		+
		K\beta\prod_{i=1}^K|w_i|^q
		\right)dx.
		\]
		
		Up to a further subsequence, we may assume that
		\[
		\int_\Omega|\nabla w_{i,n}|^2dx\to b_i\geq0,
		\qquad i=1,\ldots,K.
		\]
		Summing \eqref{eq:remainder-nehari} over $i$ and recalling that
		$Kq=2^*$, we obtain
		\[
		I(\mf w_n)
		=
		\frac1N\sum_{i=1}^K b_i+o(1).
		\]
		Since $E'(\mf u)=0$ (resp. $\mf g_{\bs s}(\mf u)=0$), we also have $E(\mf u)\geq0$. Hence,
		letting $n\to\infty$ in \eqref{eq:energy-splitting2},
		\begin{equation}\label{eq:B-splitting}
			0\leq E(\mf u)
			\leq
			E(\mf u)+\frac1N\sum_{i=1}^K b_i
			=
			\cB.
		\end{equation}
		It remains to prove that $u_i\not\equiv0$ for every $i=1,\ldots,K$.
		
		\medskip
		
		\noindent\textbf{Case 1: $\mf u=\bs0$.}
		
		\medskip
		
		In this case $w_{i,n}=u_{i,n}$ and $b_i>0$ for every $i$. Indeed,
		if $b_i=0$ for some $i$, then $u_{i,n}\to0$ in $H_0^1(\Omega)$,
		contradicting Lemma \ref{lem:uniform-bounds-components}.
		
		We claim that there exist $t_{i,n}>0$, $i=1,\ldots,K$, such that
		\[
		(t_{1,n}w_{1,n},\ldots,t_{K,n}w_{K,n})\in\mathcal N.
		\]
		For $\bs t=(t_1,\ldots,t_K)\in(0,+\infty)^K$, define
		\[
		F_{j,n}(\bs t)
		:=
		t_j^2\int_\Omega|\nabla w_{j,n}|^2dx
		-
		t_j^{2^*}\mu_j\int_\Omega|w_{j,n}|^{2^*}dx
		-
		\beta t_j^q\prod_{\ell\neq j}t_\ell^q
		\int_\Omega\prod_{\ell=1}^K|w_{\ell,n}|^qdx.
		\]
		Since $(\mf w_n)$ is bounded, after passing to a subsequence we may
		assume that, for every $j=1,\ldots,K$,
		\[
		\mu_j\int_\Omega|w_{j,n}|^{2^*}dx\to c_j>0,
		\qquad
		|\beta|\int_\Omega\prod_{\ell=1}^K|w_{\ell,n}|^qdx\to d\geq0.
		\]
		By \eqref{eq:remainder-nehari},
		\begin{equation}\label{eq:bi-ci-d}
			c_j=b_j+d.
		\end{equation}
		
		It is easy to see that we can choose $0<r<1<R$ so that, for $n$ sufficiently large and
		uniformly with respect to the remaining variables in $[r,R]$,
		\[
		F_{j,n}(\bs t)>0
		\quad\text{if }t_j=r,
		\qquad
		F_{j,n}(\bs t)<0
		\quad\text{if }t_j=R.
		\]
		Miranda's theorem therefore yields
		$\bs t_n=(t_{1,n},\ldots,t_{K,n})\in[r,R]^K$ such that
		\[
		F_{j,n}(\bs t_n)=0,
		\qquad j=1,\ldots,K.
		\]
		Up to a subsequence, $t_{i,n}\to\bar t_i\in[r,R]$.
		Passing to the limit gives
		\begin{equation}\label{eq:t-limit-system}
			\bar t_i^2 b_i
			=
			\bar t_i^{2^*}c_i
			-
			d\prod_{\ell=1}^K\bar t_\ell^q,
			\qquad i=1,\ldots,K.
		\end{equation}
		In particular,
		\begin{equation}\label{eq:t-lower-bound}
			\bar t_i
			\geq
			\left(\frac{b_i}{c_i}\right)^{\frac1{2^*-2}}.
		\end{equation}
		
		We claim that $\bar t_i=1$ for every $i$. If $d=0$, this follows
		immediately from \eqref{eq:bi-ci-d} and \eqref{eq:t-limit-system}.
		Suppose therefore that $d>0$.
		
		For every $i$, let
		\[
		f_i(s):=s^{(K-1)q}c_i-s^{2-q}b_i.
		\]
		Since $(K-1)q=2^*-q$, we have
		\[
		f_i'(s)
		=
		s^{1-q}
		\left[
		(2^*-q)c_i s^{2^*-2}-(2-q)b_i
		\right]>(2-q)s^{1-q}(c_i s^{2^*-2}-b_i)\geq0
		\]
		for
		\[
		s\geq
		\left(\frac{b_i}{c_i}\right)^{\frac1{2^*-2}}.
		\]
		Moreover, \eqref{eq:bi-ci-d} gives $f_i(1)=c_i-b_i=d$.
		Dividing \eqref{eq:t-limit-system} by $\bar t_i^q$, we obtain
		\begin{equation}\label{eq:f-system}
			d\prod_{\ell\neq i}\bar t_\ell^q
			=f_i(\bar t_i).
		\end{equation}
		
		Suppose that $\bar t_i<1$ for some $i$. Then, \eqref{eq:t-lower-bound},
		\eqref{eq:f-system} and the monotonicity of $f_i$ imply
		\[
		\prod_{\ell\neq i}\bar t_\ell<1.
		\]
		Hence at least two components are smaller than $1$. After relabelling,
		we may assume that, for some $m\geq2$,
		\[
		\bar t_1,\ldots,\bar t_m<1,
		\qquad
		\bar t_{m+1},\ldots,\bar t_K\geq1.
		\]
		Multiplying \eqref{eq:f-system} for $i=1,\ldots,m$, we obtain
		\begin{equation*}
			d^m
			=
			\left(\prod_{\ell=m+1}^K\bar t_\ell^{-mq}\right)
			\prod_{i=1}^m
			\bar t_i^{(K-m)q}
			\left(c_i-\bar t_i^{2-2^*}b_i\right)
			<
			\prod_{i=1}^m(c_i-b_i)
			=d^m,
		\end{equation*}
		a contradiction.
		If instead $\bar t_i>1$ for some $i$, then
		\eqref{eq:f-system} gives
		\[
		\prod_{\ell\neq i}\bar t_\ell>1.
		\]
		Again, after relabelling, for some $m\geq2$ we may assume
		\[
		\bar t_1,\ldots,\bar t_m>1,
		\qquad
		\bar t_{m+1},\ldots,\bar t_K\leq1.
		\]
		The same computation now gives
		\[
		d^m
		=
		\left(\prod_{\ell=m+1}^K\bar t_\ell^{-mq}\right)
		\prod_{i=1}^m
		\bar t_i^{(K-m)q}
		\left(c_i-\bar t_i^{2-2^*}b_i\right)
		>d^m,
		\]
		again a contradiction. Therefore $\bar t_i=1$ for every $i=1,\ldots,K$. Since
		\[
		(t_{1,n}w_{1,n},\ldots,t_{K,n}w_{K,n})\in\mathcal N,
		\]
		we obtain
		\[
		\cA
		\leq
		\lim_{n\to\infty}
		I(t_{1,n}w_{1,n},\ldots,t_{K,n}w_{K,n})
		=
		\lim_{n\to\infty}I(\mf w_n)
		=
		\frac1N\sum_{i=1}^K b_i.
		\]
		On the other hand, \eqref{eq:B-splitting} yields
		\[
		\frac1N\sum_{i=1}^K b_i\leq \cB,
		\]
		thus $\cB\geq \cA$, contradicting \eqref{eq:upper bound B}.
		Hence $\mf u\neq\bs0$.
		
		\medskip
		
		\noindent\textbf{Case 2: $\mf u\neq\bs0$ and at least one component of
			$\mf u$ vanishes.}
		
		\medskip
		
		After relabelling, assume that
		\[
		u_1=\cdots=u_m=0,
		\qquad
		u_i\neq0
		\quad\text{for }i=m+1,\ldots,K,
		\]
		for some $1\leq m<K$. Since one component vanishes, the interaction
		term of the limit is identically zero. Therefore, for every
		$i=m+1,\ldots,K$,
		\[
		-\Delta u_i+\lambda_i u_i
		=
		\mu_i|u_i|^{2^*-2}u_i,
		\]
		and hence
		\begin{equation}\label{eq:energy-nonzero-components}
			E(\mf u)
			\geq
			\sum_{i=m+1}^K B_{\mu_i}.
		\end{equation}
		
		For $i=1,\ldots,m$, we have $b_i>0$, by the same argument as in
		Case~1. We next claim that $b_i=0$ for at least one
		$i\in\{m+1,\ldots,K\}$. Otherwise $b_i>0$ for every $i=1,\ldots,K$,
		and the argument of Case~1 gives
		\[
		\cA
		\leq
		\frac1N\sum_{i=1}^K b_i
		\leq
		E(\mf u)+\frac1N\sum_{i=1}^K b_i
		=
		\cB,
		\]
		again contradicting \eqref{eq:upper bound B}.
		
		Thus $b_i=0$ for some $i\geq m+1$. Consequently,
		$w_{i,n}\to0$ strongly in $H_0^1(\Omega)$, and H\"older's inequality
		yields
		\[
		\int_\Omega\prod_{\ell=1}^K|w_{\ell,n}|^q\,dx\to0.
		\]
		Therefore, for every $i=1,\ldots,m$,
		\eqref{eq:remainder-nehari} gives
		\[
		\int_\Omega |\nabla w_{i,n}|^2dx
		=
		\mu_i\int_\Omega |w_{i,n}|^{2^*}dx+o(1).
		\]
		By the Sobolev inequality,
		\[
		\int_\Omega |\nabla w_{i,n}|^2dx
		\leq
		\mu_i S^{-\frac{2^*}{2}}
		\left(
		\int_\Omega|\nabla w_{i,n}|^2
		dx\right)^{\frac{2^*}{2}}
		+o(1).
		\]
		Since $b_i>0$, passing to the limit gives
		\[
		b_i
		\geq
		\mu_i^{-\frac{N-2}{2}}S^{\frac N2},
		\qquad i=1,\ldots,m.
		\]
		Combining this with \eqref{eq:B-splitting} and
		\eqref{eq:energy-nonzero-components}, we obtain
		\[
		\cB
		=
		E(\mf u)+\frac1N\sum_{i=1}^K b_i
		\geq
		\sum_{i=m+1}^K B_{\mu_i}
		+
		\frac1N\sum_{i=1}^m
		\mu_i^{-\frac{N-2}{2}}S^{\frac N2},
		\]
		contradicting \eqref{eq:upper bound B}.
		
		Therefore neither Case~1 nor Case~2 can occur, and hence
		\[
		u_i\not\equiv0
		\qquad\text{for every }i=1,\ldots,K.
		\]
		Therefore, $\mf u\in\mathcal M$. By the
		definition of $\cB$, $E(\mf u)\geq \cB$, whereas \eqref{eq:B-splitting} gives $E(\mf u)\leq \cB$. Thus
		\[
		E(\mf u)=\cB.
		\]
		Therefore, recalling Proposition \ref{prop:minimizers_are_solutions} (and Lemma \ref{lem:clarke critical point are solutions} in the case $q=1$), we deduce that $\mf u$ is a least energy fully non-trivial solution of \eqref{eq:system 1}.
		Finally, replacing each component with its absolute value, we may
		assume that $\mf u$ is nonnegative.

    \medskip

    \textit{Case $q=1$:} If, instead, $q=1$, we need to tackle the differentiability issue with more care. By Ekeland's variational principle again, there exists a sequence $(\mf u_n)_n\subset \cM$ such that
        \begin{align*}
            E(\mf u_n)&<\cB+\frac{1}{n}, \\
            E(\mf v)&\geq E(\mf u_n)-\frac{1}{n}\norm{\mf v-\mf u_n} \quad \forall \mf v\in \cM, \mf v \neq \mf u_n.
        \end{align*}
        
Therefore, setting
\[
F_n(\mf v):=
E(\mf v)+\frac1n\|\mf v-\mf u_n\|_H,
\]
we see that $\mf u_n$ minimizes $F_n$ on $\cM$.
By Clarke's generalized Lagrange multiplier rule
\cite[Theorem~1 and Corollary~1]{Clarke_1976},
there exist $\lambda_{0,n}\geq0$ and
$\lambda_{i,n}\in\R$, $i=1,\ldots,K$,
not all zero, such that
\[
0\in
\lambda_{0,n}\partial F_n(\mf u_n)
+\sum_{i=1}^K
\lambda_{i,n}\partial G_{\beta,i}(\mf u_n).
\]
If $\lambda_{0,n}=0$, testing this inclusion
against $\bs\varphi_j(\mf u_n)$,
$j=1,\ldots,K$, yields
\[
\bs M_\beta(\mf u_n)^T\bs\lambda_n=\bs0.
\]
By \Cref{lemma:M_beta_negative},
$\bs\lambda_n=\bs0$, a contradiction.
Thus $\lambda_{0,n}>0$, and after normalization
we obtain, using the subdifferential sum rule and observing that $\partial \norm{\cdot-\mf u_n}(\mf u_n)\subset B^*$, where $B^*$ is the unit ball in $H'$,
\[
0\in
\partial E(\mf u_n)
+\sum_{i=1}^K
\lambda_{i,n}\partial G_{\beta,i}(\mf u_n)
+\frac1n B^*,
\]
where $B^*$ denotes the closed unit ball in $H'$.

Recalling \eqref{eq:generalized-gradient-E}, there exist
$\bs s_n,\bs t_{i,n}\in\cS(\mf u_n)$, $i=1,\dots,K$, such that
\[
\xi_n
:=
\mf g_{\bs s_n}(\mf u_n)
-
\sum_{i=1}^K
\lambda_{i,n}
G_{\beta,i,\bs t_{i,n}}'(\mf u_n)\to 0 \quad \text{ in $H'$}.
\]
For every $j=1,\dots,K$, we test the above identity against
$\bs\varphi_j(\mf u_n)$. Since $(\mf u_n)$ is bounded in $H$, we have
\[
\xi_n[\bs\varphi_j(\mf u_n)]=o(1).
\]
Moreover, since $\mf u_n\in\mathcal M$, the definition of
$\bs\varphi_j(\mf u_n)$ gives
\[
\mf g_{\bs s_n}(\mf u_n)[\bs\varphi_j(\mf u_n)]
=
\|u_{j,n}\|_{\lambda_j}^2
-\mu_j|u_{j,n}|_{2^*}^{2^*}
-\beta\int_\Omega\prod_{k=1}^K|u_{k,n}|^q\,dx
=
G_{\beta,j}(\mf u_n)
=
0.
\]
Here, in the case $q=1$, we have used that $s_{j,n}u_{j,n}=|u_{j,n}|$ a.e. in $\Omega$, so that the value of the interaction term is independent of the
particular choice of $\bs s_n\in\cS(\mf u_n)$.

Similarly, by the definition of $\bs M_\beta(\mf u_n)$ and the
independence of its entries from the particular choice of
$\bs t_{i,n}\in\cS(\mf u_n)$ when $q=1$, we have
\[
G_{\beta,i,\bs t_{i,n}}'(\mf u_n)
[\bs\varphi_j(\mf u_n)]
=
(\bs M_\beta(\mf u_n))_{ij}.
\]
Therefore,
\[
\sum_{i=1}^K
\lambda_{i,n}
(\bs M_\beta(\mf u_n))_{ij}
=o(1),
\qquad j=1,\dots,K,
\]
or, equivalently,
\[
\bs M_\beta(\mf u_n)^T\bs\lambda_n=o(1)
\qquad\text{in }\R^K,
\]
where
\[
\bs\lambda_n
:=
(\lambda_{1,n},\dots,\lambda_{K,n})^T.
\]
By the uniform invertibility of $\bs M_\beta(\mf u_n)$ established
in Lemma \ref{lemma:M_beta_negative}, we conclude that
\[
\bs\lambda_n\to\bs0.
\]
Since the generalized gradients
$\partial G_{\beta,i}(\mf u_n)$ are uniformly bounded in $H'$ along
bounded subsets of $H$, it follows that
\[
\dist(0,\partial E(\mf u_n))\to0.
\]
Therefore, there exists
$\bs s_n=(s_{1,n},\dots,s_{K,n})\in\cS(\mf u_n)$ such that $\xi_n=\mf g_{\bs s_n}(\mf u_n)\to 0$ in $H'$ as $n\to\infty$.
Hence, for every
$\bs\varphi=(\varphi_1,\dots,\varphi_K)\in C_c^\infty(\Omega,\R^K)$,
\begin{align}
o(1)
&=
\xi_n[\bs\varphi]
\nonumber\\
&=
\sum_{i=1}^K
\int_\Omega
\left(
\nabla u_{i,n}\cdot\nabla\varphi_i
+\lambda_i u_{i,n}\varphi_i
\right)\,dx
-
\sum_{i=1}^K
\mu_i
\int_\Omega
|u_{i,n}|^{2^*-2}u_{i,n}\varphi_i\,dx
\nonumber\\
&\quad
-
\beta\sum_{i=1}^K
\int_\Omega
s_{i,n}
\prod_{j\neq i}|u_{j,n}|
\,\varphi_i\,dx .
\label{eq:PS-generalized-limit}
\end{align}
By the weak convergence in $H_0^1(\Omega)$, the compact Sobolev
embeddings, and the a.e.\ convergence, all the terms in
\eqref{eq:PS-generalized-limit} except for the interaction term pass
immediately to the limit.

We therefore consider the latter. Since $|s_{i,n}|\leq1$ a.e.\ in
$\Omega$, up to a subsequence,
\[
s_{i,n}\overset{*}{\rightharpoonup}s_i
\qquad\text{in }L^\infty(\Omega)
\]
for some $s_i\in L^\infty(\Omega)$ with $|s_i|\leq1$ a.e.. Moreover,
\[
s_i(x)\in\operatorname{Sign}(u_i(x))
\qquad\text{for a.e. }x\in\Omega.
\]
Indeed, since $u_{i,n}\to u_i$ a.e.\ in $\Omega$, we have
$s_{i,n}\to1$ a.e.\ on $\{u_i>0\}$ and
$s_{i,n}\to-1$ a.e.\ on $\{u_i<0\}$. Hence, by dominated convergence
and the weak-$*$ convergence of $s_{i,n}$ in $L^\infty(\Omega)$,
\[
s_i=1\quad\text{a.e. on }\{u_i>0\},
\qquad
s_i=-1\quad\text{a.e. on }\{u_i<0\}.
\]
Moreover, $|s_i|\leq1$ a.e.\ in $\Omega$, and therefore
$s_i\in\operatorname{Sign}(u_i)$ a.e.\ in $\Omega$. Thus
$\bs s=(s_1,\dots,s_K)\in\cS(\mf u)$.

Since $q=1$, we have $K=2^*$. By the compact Sobolev embedding,
\[
u_{j,n}\to u_j
\qquad\text{strongly in }L^{K-1}(\Omega),
\]
and therefore
\[
\prod_{j\neq i}|u_{j,n}|
\longrightarrow
\prod_{j\neq i}|u_j|
\qquad\text{strongly in }L^1(\Omega).
\]
Since $\varphi_i\in C_c^\infty(\Omega)$, it follows that
\begin{align*}
&\int_\Omega
s_{i,n}\prod_{j\neq i}|u_{j,n}|\varphi_i\,dx
-
\int_\Omega
s_i\prod_{j\neq i}|u_j|\varphi_i\,dx
\\
&\quad=
\int_\Omega
s_{i,n}
\left(
\prod_{j\neq i}|u_{j,n}|
-
\prod_{j\neq i}|u_j|
\right)\varphi_i\,dx
+
\int_\Omega
(s_{i,n}-s_i)
\prod_{j\neq i}|u_j|\varphi_i\,dx
\longrightarrow0,
\end{align*}
by the strong $L^1$ convergence and the weak-$*$ convergence in
$L^\infty(\Omega)$, respectively.

Passing to the limit in \eqref{eq:PS-generalized-limit}, we obtain
\[
\mf g_{\bs s}(\mf u)[\bs\varphi]=0
\qquad
\text{for every }\bs\varphi\in C_c^\infty(\Omega,\R^K).
\]
By density and continuity, the same identity holds for every
$\bs\varphi\in H$. Therefore
\[
\mf g_{\bs s}(\mf u)=0
\qquad\text{for some }\bs s\in\cS(\mf u).
\]
From now on, the proof proceeds as in the previous case.
\end{proof}

\section{Examples of domains satisfying Condition \eqref{eq:condition-C}}\label{sec:domains condition C}
We now show that Condition \eqref{eq:condition-C} can be satisfied by constructing a class of elongated domains $\Omega_{L,m} \subset \mathbb{R}^N$ ($N=3,4$) for which \eqref{eq:condition-C} holds when $L > 0$ is sufficiently large, for every admissible number of components $K\geq 3$. We also show that condition \eqref{eq:lambda-general} is satisfied, since this is an essential ingredient in the existence theorems.

For $L > 1$ and $m \in \mathbb{N}$, $m \ge 1$, we define the domain
\begin{equation*}
\Omega_{L,m} := \left\{ (x_1, x') \in \mathbb{R} \times \mathbb{R}^{N-1} : \left( \frac{x_1}{L} \right)^{2m} + |x'|^2 < 1 \right\}.
\end{equation*}
Let $\omega_{N-1} := |B_1^{N-1}|$ denote the volume of the unit ball in $\mathbb{R}^{N-1}$, and let $\nu_N := \lambda_1(B_1^{N-1})$ be its first Dirichlet eigenvalue.

\begin{lemma}\label{lem:domain-properties}
For every $L > 1$ and $m \in \mathbb{N}$, the domain $\Omega_{L,m}$ satisfies:
\begin{enumerate}
    \item \textbf{Volume estimate:} $|\Omega_{L,m}| \leq 2 L \omega_{N-1}$.
    \item \textbf{First eigenvalue:} $\Lambda_1(\Omega_{L,m}) \geq \nu_N + \frac{\pi^2}{4L^2}$.
    \item \textbf{Second eigenvalue:} for every $s \in (0,1)$, we have
    \[
    \Lambda_2(\Omega_{L,m}) \leq \frac{\nu_N}{1 - s^{2m}} + \frac{\pi^2}{s^2 L^2}.
    \]
\end{enumerate}
\end{lemma}

\begin{proof}
Estimate (1) follows directly from the inclusion $\Omega_{L,m} \subset C_L := (-L, L) \times B_1^{N-1}$. Property (2) follows from the monotonicity of eigenvalues with respect to domain inclusion, using separation of variables in $C_L$:
\[
\Lambda_1(\Omega_{L,m}) \geq \Lambda_1(C_L) = \lambda_1((-L,L)) + \lambda_1(B_1^{N-1}) = \frac{\pi^2}{4L^2} + \nu_N.
\]
For property (3), we set $r_s := \sqrt{1 - s^{2m}}$, and
\[
Q_s =(-sL,sL) \times B_{r_s}^{N-1}.
\]
It is immediate that $Q_s \subset \Omega_{L,m}$ for every $s$. Thus,
by domain monotonicity and scaling of eigenvalues,
\[
\begin{split}
\Lambda_2(\Omega_{L,m}) &\leq \Lambda_2\left(Q_s\right)  = \min\left\{\Lambda_1\left(-sL,sL\right) + \Lambda_2( B_{r_s}^{N-1}), \ \Lambda_2\left(-sL,sL\right) + \Lambda_1( B_{r_s}^{N-1})\right\}\\
& \le \Lambda_2\left(-sL,sL\right) + \Lambda_1( B_{r_s}^{N-1}) = \frac{\pi^2}{s^2 L^2} + \frac{\nu_N}{1-s^{2m}},
\end{split}
\]
which proves (3).
\end{proof}

\begin{proposition}\label{prop:condition-C-satisfied}
Let $N \in \{3,4\}$ and $m \in \N$, $m \ge 1$. Assume that $\lambda_i = -\nu_N$ for every $i$.
Then, there exists $L_0>0$ such that Condition \eqref{eq:condition-C} holds on $\Omega_{L,m}$ for all $L>L_0$. In addition, with this choice also assumption \eqref{eq:lambda-general} is satisfied.
\end{proposition}

\begin{proof}
    Using Lemma \ref{lem:domain-properties} and the fact that $\overline{\lambda}=-\nu_N$, we have the estimate
\[
\Lambda_2(\Omega_{L,m}) + \overline{\lambda} \leq \frac{s^{2m}\nu_N}{1-s^{2m}} + \frac{\pi^2}{s^2 L^2},
\]
for every $s \in (0,1)$. Now we choose $s=L^{-\frac1{m+1}}$, so that $s^2L^2 = s^{-2m}$, and by taking $L$ sufficiently large
\[
\Lambda_2(\Omega_{L,m}) + \overline{\lambda} \le (2\nu_N + \pi^2)s^{2m}.
\]
Note also that the left hand side is positive, since $\Lambda_2(\Omega_{L,m}) > \Lambda_1(\Omega_{L,m}) > \nu_N$. With the preceding estimate at hand, we can estimate the left hand side in Condition \eqref{eq:condition-C} as follows:
\beq\label{check cond C}
\frac{1}{N} \overline{\mu} \left( \Lambda_2(\Omega_{L,m}) + \overline{\lambda} \right)^{N/2} |\Omega_{L,m}| \leq C L^{1 - \frac{Nm}{m+1}} < \frac{1}{N} \underline{\mu} S^{N/2}
\eeq
for all sufficiently large $L>L_0$, where $C$ denotes a positive constant $C>0$ depending only on $\bar \mu$ and the dimension, and we used the fact that the exponent of $L$ in the second term is negative, since $m+1>0$ and $m+1<Nm$. This proves that Condition \eqref{eq:condition-C} is fulfilled for any such $L$.

Regarding assumption \eqref{eq:lambda-general}, by point (2) of
Lemma \ref{lem:domain-properties} we have that
$-\Lambda_1(\Omega_{L,m})<\lambda_i<0$ for every $i$, both in dimension $3$
and $4$. Thus, it remains to check that
$\lambda_i=-\nu_N<-\lambda^*(\Omega_{L,m})$ when $N=3$. Recall that
\[
\lambda^*(\Omega_{L,m})=\inf\{\lambda>0:\ S_\lambda(\Omega_{L,m})<S\},
\]
where, for $\lambda<\Lambda_1(\Omega_{L,m})$,
\[
S_\lambda(\Omega_{L,m}):=\inf_{u\in H^1_0(\Omega_{L,m})\setminus\{0\}}
\frac{\|\nabla u\|_2^2-\lambda\|u\|_2^2}{\|u\|_{2^*}^2}
\]
(so that the least energy level of the Brezis-Nirenberg problem
$-\Delta u=\lambda u+|u|^{2^*-2}u$ in $\Omega_{L,m}$ is
$S_\lambda(\Omega_{L,m})^{N/2}/N$). Testing the definition of
$S_\lambda(\Omega_{L,m})$ with the first eigenfunction $\varphi_1$, we obtain
\begin{equation}\label{eq:test-eigenfunction}
S_\lambda(\Omega_{L,m})\le\big(\Lambda_1(\Omega_{L,m})-\lambda\big)\,
|\Omega_{L,m}|^{\frac2N}
\qquad\text{for every }\lambda<\Lambda_1(\Omega_{L,m}).
\end{equation}
On the other hand, since $\overline\lambda=-\nu_N$,
$0<\Lambda_1(\Omega_{L,m})-\nu_N\le\Lambda_2(\Omega_{L,m})-\nu_N$
(by point (2) of Lemma \ref{lem:domain-properties}) and
$\underline\mu/\overline\mu\le1$, estimate \eqref{check cond C} gives
\begin{equation}\label{eq:gap-volume}
\big(\Lambda_1(\Omega_{L,m})-\nu_N\big)\,|\Omega_{L,m}|^{\frac2N}<S
\qquad\forall L>L_0 .
\end{equation}
Taking $\lambda=\nu_N$ in \eqref{eq:test-eigenfunction}, this already yields
$S_{\nu_N}(\Omega_{L,m})<S$; however, since $\lambda^*(\Omega_{L,m})$ is
defined as an infimum, we need a little room to get a strict inequality.
Fix $L>L_0$ and set
\[
\varepsilon:=\frac12\min\Big\{\nu_N,\ S\,|\Omega_{L,m}|^{-\frac2N}
-\big(\Lambda_1(\Omega_{L,m})-\nu_N\big)\Big\},
\]
which is positive by \eqref{eq:gap-volume}. Then
$\lambda':=\nu_N-\varepsilon\in(0,\nu_N)$, and in particular
$\lambda'<\Lambda_1(\Omega_{L,m})$. By \eqref{eq:test-eigenfunction} and the
choice of $\varepsilon$,
\[
S_{\lambda'}(\Omega_{L,m})
\le\big(\Lambda_1(\Omega_{L,m})-\nu_N+\varepsilon\big)\,|\Omega_{L,m}|^{\frac2N}
\le\frac12\Big(S+\big(\Lambda_1(\Omega_{L,m})-\nu_N\big)|\Omega_{L,m}|^{\frac2N}\Big)
<S .
\]
By the definition of $\lambda^*(\Omega_{L,m})$ as an infimum, we infer
$\lambda^*(\Omega_{L,m})\le\lambda'<\nu_N$, that is,
\[
\lambda_i=-\nu_N<-\lambda^*(\Omega_{L,m})\qquad\text{for every }i .
\]
Therefore, \eqref{eq:lambda-general} also holds for $N=3$, and the proof is
complete.
\end{proof}

\begin{remark}
    It is not necessary that all the $\lambda_i$'s are equal. Since \eqref{eq:condition-C} is an open condition, for $m$ and $L$ such that Proposition \ref{prop:condition-C-satisfied} holds we can take $\lambda_i =-\nu_N + \eta_i$ with $|\eta_i|$ sufficiently small.
\end{remark}

\section{Asymptotic behavior under strong competition}\label{sec:asymptotic behaviour}

In this section, we investigate the limiting behavior, as
$\beta\to-\infty$, of the nonnegative least energy fully non-trivial
solutions obtained in
Theorems~\ref{thm:least energy competitive N=5,6}
and~\ref{thm:least energy competitive N=3,4},
and prove Theorem~\ref{thm:limit behaviour}.
To emphasize the dependence on $\beta$, throughout this section
we write $E_\beta$, $\cB_\beta$, and $\cM_\beta$
in place of $E$, $\cB$, and $\cM$, respectively.

We begin with the following lemma, which strengthens
\Cref{lem:upper-bound-B-N56} in dimension $N=6$.

\begin{lemma}\label{lem:bound N=6}
		Let $K=3$, $\beta<0$, and $N=6$. Then
		\begin{equation*}
			\sup_{\b<0} \cB_\beta<
			\min\left\{
		\min_{1\leq \ell\leq 3}\Big(\sum_{k\not=\ell} B_{\mu_k}
		+\frac1N\mu_\ell^{-\frac{N-2}{2}}S^{\frac N2}\Big),\cA
		\right\}.
		\end{equation*}
	\end{lemma}

    \begin{proof}
 By  \eqref{eq:A value beta negative} and \eqref{eq:B_mu less than S}, it is enough to prove
    $$
    \sup_{\b<0} \cB_\beta<B_{\mu_1}+B_{\mu_2}+\frac1N\mu_3^{-\frac{N-2}{2}}S^{\frac N2}.
    $$
    Since the  main argument follows from \cite[Lemma 5.1]{Chen_Zou_2015}, we only sketch the proof here.
    Indeed, for the problem
    \begin{equation}\label{eq:B_R} 
\begin{cases}
-\Delta u+\lambda_3 u
=\mu_3|u|^{2^*-2}u
&\text{in }B(0,R),\\
u>0,~u\in H_0^1(B(0,R)),
\end{cases} 
\end{equation}
we define the associated energy functional $J_R:H_0^1(B(0,R))\rightarrow\mathbb{R}$ as
\[
J_R(u)=\frac{1}{2}\int_{B(0,R)}(|\nabla u|^2+\lambda_3 u^2)dx-\frac{1}{2^*}\mu_3\int_{B(0,R)}|u|^{2^*}dx.
\]
It follows from \cite[Theorem 5.1]{Chen_Zou_2015} that there is $R_0>0$ such that for any $R\in(0,R_0)$,
\eqref{eq:B_R} admits a least energy solution $U_R$ and
\[
\frac{1}{N}\mu_3^{-\frac{N-2}{2}}S^{N/2}-C_1R^{\frac{2N-4}{N-4}}\leq J_R(U_R)\leq \frac{1}{N}\mu_3^{-\frac{N-2}{2}}S^{N/2}-C_2R^{\frac{2N-4}{N-4}},
\]
for some $C_1, C_2>0$. Recall that $u_{\mu_i}$ denotes a positive least energy solution of the scalar Brezis-Nirenberg problem \eqref{Brezis-Nirenberg} for $i=1,2,3$, let $R>0$ sufficiently small and $x_R\in \Omega$ such that
$\dist(x_R,\pa \Omega)=3R$ and 
\[
|u_{\mu_1}(x)|\leq CR,~~ x\in B(x_R,3R).
\]
Denote also $\varphi_R:=1-\psi(\frac{x-x_R}{R})$, where $\psi\in C_0^\infty(B(0,2))$ with $0\leq \psi\leq 1$ and $\psi\equiv1$ in $B(0,1)$.
Set $u_{1,R}:=\varphi_R u_{\mu_1}$, $u_{2,R}:=u_{\mu_2}$, and $u_{3,R}:=U_R(\cdot-x_R)$. Then,
similarly as in  \cite[Lemma 5.1]{Chen_Zou_2015}, we can find $t_R>0$ such that
    $$
    (t_R u_{1,R}, u_{2,R}, u_{3,R})\in \cM_\b
    $$
    for every $\b<0$, and estimate the energy of this competitor to obtain the desired result.
    \end{proof}

    Finally, we shall prove Theorem \ref{thm:limit behaviour}.
    \begin{proof}[Proof of Theorem \ref{thm:limit behaviour}.]
  Let  $\mf u_n:=(u^{(n)}_1,\dots,u^{(n)}_K)\in H$ be the nonnegative solutions of \eqref{eq:system 1} with $\beta=\beta_n$ given by Theorem \ref{thm:least energy competitive N=5,6} or Theorem \ref{thm:least energy competitive N=3,4}, where $\beta_n\rightarrow-\infty$ as $n\rightarrow\infty$. Then, it follows from Lemma \ref{lem:upper-bound-B-N56} when \(N=5,6\), and Lemma \ref{lem:abstract-bound-N34} when \(N=3,4\), that
\begin{equation*}
\cB_{\beta_n}=E_{\beta_n}(\mf u_n)<\displaystyle\sum_{i=1}^{K}\frac1N\mu_i^{-\frac{N-2}{2}}S^{\frac N2},
\end{equation*}
which immediately implies that $u^{(n)}_i$ is uniformly bounded in $H_0^1(\Omega)$ with respect to $n$, where $i=1,\cdots,K$. Up to a subsequence, if necessary, we obtain $u^{(n)}_i\rightharpoonup u^{\infty}_i$ weakly in $H_0^1(\Omega)$, $u^{(n)}_i\rightarrow u^{\infty}_i$ strongly in $L^2(\Omega)$ and $u^{(n)}_i\rightarrow u^{\infty}_i\geq0$ a.e. for $i=1,\cdots,K$.
Let us show, in three steps, that either conclusion (1) or conclusion (2) hold.
\medskip

\noindent\textbf{Step 1:} it is not possible that $u^{\infty}_i\equiv0$ for all $i\in\{1,\cdots,K\}$. 

In that case, $u^{(n)}_i\to 0$ in $L^2(\Omega)$ for all $i$. Noticing that $\mf u_n\in \cM_{\beta_n}$ and $\beta_n<0$, it is easy to see that
\begin{equation*}
\displaystyle\int_{\Omega}(|\nabla u^{(n)}_i|^2+\lambda_i|u^{(n)}_i|^2)dx\leq
\displaystyle\int_{\Omega}\mu_i|u^{(n)}_i|^{2^{\ast}}dx\leq \mu_iS^{-\frac{2^{\ast}}{2}}\left(\displaystyle\int_{\Omega}|\nabla u^{(n)}_i|^2dx\right)^{\frac{2^{\ast}}{2}}
\end{equation*}
and so
\begin{equation}\label{eq6.3}
b_i:=\displaystyle\lim_{n\rightarrow\infty}\displaystyle\int_{\Omega}|\nabla u^{(n)}_i|^2dx\geq\mu_i^{-\frac{N-2}{2}}S^{\frac N2}
\end{equation}
for every $i$. In view of Lemmas \ref{lem:upper-bound-B-N56}, \ref{lem:abstract-bound-N34} and (\ref{eq6.3}), we have
\begin{align*}
\displaystyle\sum_{i=1}^{K-1}B_{\mu_i}+\frac{1}{N}\mu_K^{-\frac{N-2}{2}}S^{\frac N2}&\geq\displaystyle\lim_{n\rightarrow\infty}\cB_{\beta_n}=\frac{1}{N}\displaystyle\sum_{i=1}^{K}b_i
\\&\geq\frac{1}{N}\displaystyle\sum_{j=1}^{K}\mu_j^{-\frac{N-2}{2}}S^{\frac N2}>\displaystyle\sum_{i=1}^{K-1}B_{\mu_i}+\frac{1}{N}\mu_K^{-\frac{N-2}{2}}S^{\frac N2},
\end{align*}
which is a contradiction.

\medskip

\noindent\textbf{Step 2:} if $u^{\infty}_i\equiv0$ for at least one index $i_0$, then necessarily
$$
i_0\in \argmin_{j=1,\dots,K} \left(\displaystyle\sum_{i\neq j}B_{\mu_i}+\frac{1}{N}\mu_j^{-\frac{N-2}{2}}S^{\frac N2}\right)
$$
and $u_j^\infty\not\equiv 0$ for every $j\neq i_0$. In that case, conclusion $(1)$ and \eqref{eq:vanishing interacion} hold.

\medskip

We set $I_0:=\{i\in \{1,\dots,K\}: u_i^\infty\equiv 0\}$. By Step 1, at least one component $u_i^\infty\not\equiv 0$, therefore $\{1,\dots,K\}\setminus I_0\neq \emptyset$. First, arguing as before,
\begin{equation}\label{eq:b_i}
b_i:=\displaystyle\lim_{n\rightarrow\infty}\displaystyle\int_{\Omega}|\nabla u^{(n)}_i|^2dx\geq\mu_i^{-\frac{N-2}{2}}S^{\frac N2}
\end{equation}
for every $i\in I_0$. Recalling that $\mf u_n$ is a solution of \eqref{eq:system 1} and $\beta_n<0$, we have
  \begin{equation*}
\displaystyle\int_{\Omega}\nabla u_i^{(n)}\nabla u_i^{\infty}dx+\lambda_i\displaystyle\int_{\Omega} u_i^{(n)}u_i^{\infty}dx\leq\mu_i\displaystyle\int_{\Omega}(u_i^{(n)})^{2^{\ast}-1}u_i^{\infty}dx
\end{equation*}
for $i\not\in I_0$. Passing to the limit using the convergences above and recalling the definition of $B_{\mu_i}$ shows that
\begin{equation}\label{eq6.4}
\displaystyle\int_{\Omega}(|\nabla u^{\infty}_i|^2+\lambda_i|u^{\infty}_i|^2)dx\leq
\displaystyle\int_{\Omega}\mu_i|u^{\infty}_i|^{2^{\ast}}dx\leq (NB_{\mu_i})^{-\frac{2^{\ast}}{N}}\left(\displaystyle\int_{\Omega}(|\nabla u^{\infty}_i|^2+\lambda_i|u^{\infty}_i|^2)dx\right)^{\frac{2^{\ast}}{2}},
\end{equation}
and so
\begin{equation}\label{eq6.5}
\displaystyle\int_{\Omega}(|\nabla u^{\infty}_i|^2+\lambda_i|u^{\infty}_i|^2)dx\geq NB_{\mu_i}
\end{equation}
for $i\not \in I_0$. 

By Lemmas \ref{lem:upper-bound-B-N56} and \ref{lem:abstract-bound-N34}, and recalling that $B_{\mu_i}<\frac{1}{N}\mu_i^{-\frac{N-2}{2}}S^{\frac N2}$ for every $i$, it follows that
\begin{equation}\label{eq:chain}
    \begin{aligned}
        \min_{1\leq j\leq K} &\left(\displaystyle\sum_{i\neq j}B_{\mu_i}+\frac{1}{N}\mu_j^{-\frac{N-2}{2}}S^{\frac N2}\right)\geq\displaystyle\lim_{n\rightarrow\infty}\cB_{\beta_n}=\displaystyle\lim_{n\rightarrow\infty}\frac{1}{N}\sum_{i=1}^{K}\displaystyle\int_{\Omega}(|\nabla u^{(n)}_i|^2+\lambda_i|u^{(n)}_i|^2)dx
        \\&=\frac{1}{N}\sum_{i\not\in I_0}\displaystyle\int_{\Omega}(|\nabla u^{\infty}_i|^2+\lambda_i|u^{\infty}_i|^2)dx
        +\displaystyle\lim_{n\rightarrow\infty}\frac{1}{N}\sum_{i\not\in I_0}\displaystyle\int_{\Omega}|\nabla (u^{(n)}_i-u^{\infty}_i)|^2dx+\frac{1}{N}\sum_{i\in I_0}^{K}b_i
        \\&\geq\sum_{i\not\in I_0}B_{\mu_i}+\displaystyle\lim_{n\rightarrow\infty}\frac{1}{N}\sum_{i\not\in I_0}\displaystyle\int_{\Omega}|\nabla (u^{(n)}_i-u^{\infty}_i)|^2dx+\frac{1}{N}\sum_{i\in I_0}\mu_i^{-\frac{N-2}{2}}S^{\frac N2} \\
        &\geq \min_{1\leq j\leq K} \left(\displaystyle\sum_{i\neq j}B_{\mu_i}+\frac{1}{N}\mu_j^{-\frac{N-2}{2}}S^{\frac N2}\right)+\lim_{n\rightarrow\infty}\frac{1}{N}\sum_{i\not\in I_0}\displaystyle\int_{\Omega}|\nabla (u^{(n)}_i-u^{\infty}_i)|^2dx.
    \end{aligned}
\end{equation}
This implies 
$$
\lim_{n\to\infty} \displaystyle\int_{\Omega}|\nabla (u^{(n)}_i-u^{\infty}_i)|^2dx = 0 \quad \text{for every $i\not\in I_0$},
$$
i.e., $u^{(n)}_i\rightarrow u^{\infty}_i$ strongly in $H_0^1(\Omega)$. Moreover, all the preceding inequalities are equalities, and this implies
$$
I_0=\{i_0\} \quad \text{ with } i_0\in \argmin_{j=1,\dots,K}\left(\displaystyle\sum_{i\neq j}B_{\mu_j}+\frac{1}{N}\mu_j^{-\frac{N-2}{2}}S^{\frac N2}\right).
$$
Moreover, the equality implies equality in \eqref{eq6.4} and \eqref{eq6.5}, whence
$$
\displaystyle\int_{\Omega}(|\nabla u^{\infty}_i|^2+\lambda_i|u^{\infty}_i|^2)dx=\mu_i\displaystyle\int_{\Omega}|u^{\infty}_i|^{2^{\ast}}dx=NB_{\mu_i} \quad \text{for every } i\neq i_0.
$$
Thus, $u^{\infty}_i$ is a least energy solution of
\begin{equation*}
-\Delta u+\lambda_i u = \mu_i u^{2^\ast-1},~u\in H_0^1(\Omega) \quad \text{for every } i\neq i_0,
\end{equation*}
and the maximum principle shows that $u^{\infty}_i>0$ in $\Omega$
for every $i\neq i_0$. Thus $(1)$ of Theorem \ref{thm:limit behaviour} holds.

We are left to show that \eqref{eq:vanishing interacion} holds in this case. Suppose by contradiction that it is not true, since $\mf u_n$ is bounded in $H$ and $\beta_n<0$, we may assume that $$
\displaystyle\lim_{n\rightarrow\infty}\displaystyle\int_{\Omega}\beta_n\prod_{i=1}^{K}|u^{(n)}_i|^qdx=-C<0.
$$
Then
\begin{align*}
\displaystyle\int_{\Omega}(|\nabla u^{(n)}_{i_0}|^2+\lambda_{i_0}|u^{(n)}_{i_0}|^2)dx&=\displaystyle\int_{\Omega}\mu_{i_0}|u^{(n)}_{i_0}|^{2^{\ast}}dx+\displaystyle\int_{\Omega}\beta_n\prod_{i=1}^{K}|u^{(n)}_i|^qdx
\\&\leq\mu_{i_0}S^{-\frac{2^{\ast}}{2}}\left(\displaystyle\int_{\Omega}|\nabla u^{(n)}_{i_0}|^2dx\right)^{\frac{2^{\ast}}{2}}+\displaystyle\int_{\Omega}\beta_n\prod_{i=1}^{K}|u^{(n)}_i|^qdx,
\end{align*}
and so
\begin{equation}\label{eq6.8}
b_{i_0}\leq\mu_{i_0}S^{-\frac{2^{\ast}}{2}}b_{i_0}^{\frac{2^{\ast}}{2}}-C<\mu_{i_0}S^{-\frac{2^{\ast}}{2}}b_{i_0}^{\frac{2^{\ast}}{2}}.
\end{equation}
However, equality in \eqref{eq:chain} implies equality in \eqref{eq:b_i} for $i=i_0$, therefore $b_{i_0}=\mu_{i_0}^{-\frac{N-2}{2}}S^{\frac N2}$, which contradicts \eqref{eq6.8}. Thus, \eqref{eq:vanishing interacion} holds. 

\medskip

\noindent\textbf{Step 3:} if $u^{\infty}_i\not\equiv0$, for all $i\in\{1,\cdots,K\}$, then conclusion $(2)$ and \eqref{eq:vanishing interacion} hold.

\medskip

Noticing that $u^{(n)}_i$ is bounded in  $H_0^1(\Omega)$ for all $i\in\{1,\cdots,K\}$ and $\beta_n<0$, it is easy to see that  
$$
|\beta_n|\displaystyle\int_{\Omega}\prod_{i=1}^{K}|u^{(n)}_i|^qdx\leq\displaystyle\int_{\Omega}\mu_i|u^{(n)}_i|^{2^{\ast}}dx\leq C,
$$
for some $C>0$. Applying the Fatou's Lemma, we obtain
\begin{equation*}
\displaystyle\int_{\Omega}\prod_{i=1}^{K}|u^{\infty}_i|^qdx\leq\displaystyle\liminf_{n\rightarrow\infty}\displaystyle\int_{\Omega}\prod_{i=1}^{K}|u^{(n)}_i|^qdx=0.
\end{equation*}
Thus 
\begin{equation*}
E_{\beta_n}(t_1u^{\infty}_1,\cdots,t_Ku^{\infty}_K)
=\displaystyle\sum_{i=1}^{K}\left(\frac{t_i^2}{2}\displaystyle\int_{\Omega}(|\nabla u^{\infty}_i|^2+\lambda_i|u^{\infty}_i|^2)dx-\frac{t_i^{2^{\ast}}}{2^{\ast}}\displaystyle\int_{\Omega}\mu_i|u^{\infty}_i|^{2^{\ast}}dx\right):=\displaystyle\sum_{i=1}^{K} f_i(t_i).
\end{equation*}
Direct computation yields that
\begin{equation*}
\displaystyle\max_{t_i\geq0}f_i(t_i)=\frac{1}{N}(\hat{t}_i)^2\displaystyle\int_{\Omega}(|\nabla u^{\infty}_i|^2+\lambda_i|u^{\infty}_i|^2)dx=f_i(\hat{t}_i),
\end{equation*}
where
\begin{equation*}
\hat{t}_i:=\left(\frac{\displaystyle\int_{\Omega}(|\nabla u^{\infty}_i|^2+\lambda_i|u^{\infty}_i|^2)dx}{\displaystyle\int_{\Omega}\mu_i|u^{\infty}_i|^{2^{\ast}}dx}\right)^\frac{1}{2^{\ast}-2},~i=1,\cdots,K.
\end{equation*}
Arguing as above we can obtain \eqref{eq6.4}, therefore $\hat{t}_i\leq 1$ for every $i=1,\cdots,K$.
By construction, $(\hat{t}_1u^{\infty}_1,\cdots,\hat{t}_Ku^{\infty}_K)\in \cM_{\beta_n}$ for every $n$,
hence 
$$\cB_{\beta_n}\leq \displaystyle\sum_{i=1}^{K}f_i(\hat{t}_i).
$$
This shows that
\begin{align}
\displaystyle\lim_{n\rightarrow\infty}\cB_{\beta_n}&\leq \displaystyle\sum_{i=1}^{K}f_i(\hat{t}_i)=\displaystyle\sum_{i=1}^{K}\frac{(\hat{t}_i)^2}{N}\displaystyle\int_{\Omega}(|\nabla u^{\infty}_i|^2+\lambda_i|u^{\infty}_i|^2)dx\notag
\\&\leq\displaystyle\sum_{i=1}^{K}\frac{1}{N}\displaystyle\int_{\Omega}(|\nabla u^{\infty}_i|^2+\lambda_i|u^{\infty}_i|^2)dx\notag
\\&\leq\displaystyle\sum_{i=1}^{K}\frac{1}{N}\displaystyle\int_{\Omega}(|\nabla u^{\infty}_i|^2+\lambda_i|u^{\infty}_i|^2)dx+\frac{1}{N}\displaystyle\lim_{n\rightarrow\infty}\sum_{i=1}^{K}\displaystyle\int_{\Omega}|\nabla(u^{(n)}_i-u^{\infty}_i)|^2dx\notag
\\&=\displaystyle\lim_{n\rightarrow\infty}\sum_{i=1}^{K}\frac{1}{N}\displaystyle\int_{\Omega}(|\nabla u^{(n)}_i|^2+\lambda_i|u^{(n)}_i|^2)dx=\displaystyle\lim_{n\rightarrow\infty}\cB_{\beta_n}.\label{eq6.11}
\end{align}
As a result, $\hat{t}_i=1$ and $u^{(n)}_i\rightarrow u^{\infty}_i$ strongly in  $H_0^1(\Omega)$ for every $i$. Furthermore, since $(u^{\infty}_1,\cdots,u^{\infty}_K)\in \cM_{\beta_n}$, we have
\begin{equation*}
\displaystyle\lim_{n\rightarrow\infty}\displaystyle\int_{\Omega}\beta_n\prod_{i=1}^{K}|u^{(n)}_i|^qdx=\lim_{n\to+\infty} \left(\int_{\Omega}(|\nabla u^{(n)}_i|^2+\lambda_i|u^{(n)}_i|^2)dx-\int_{\Omega}\mu_i|u^{(n)}_i|^{2^{\ast}}dx\right)=0.
\end{equation*}

Next, we show that $\mf{u}_\infty:=(u^{\infty}_{1},\dots,u^{\infty}_{K})$ is a minimizer of $E_\infty$ in $\mathcal{M}_\infty$. Clearly, $\mathcal{M}_\infty\subset \mathcal{M}_{\beta_n}$ and $E_{\beta_n}(\mf{v})=E_\infty(\mf{v})$ for every $\mf v\in \cM_\infty$, thus
\[
\mathcal{B}_{\beta_n}\leq \inf_{\cM_\infty} E_\infty.
\]
On the other hand, recalling \eqref{eq6.11}, 
\[
E_\infty(\mf{u}_\infty)=\lim_{n\to\infty} \cB_{\beta_n}.
\]
This proves that $\mf{u}_\infty$ is a minimizer of
$E_\infty$ in $\mathcal{M}_\infty$.

Finally, by Lemma \ref{lem:abstract-bound-N34} when \(N=3,4\), and by Lemma \ref{lem:bound N=6} when \(N=6\),
we have
\[
\sup_{\beta<0}\mathcal{B}_\beta<\min_{1\leq \ell\leq K}
\left(
\sum_{i\ne\ell}B_{\mu_i}
+\frac1N\mu_\ell^{-\frac{N-2}{2}}S^{N/2}
\right),
\]
and this yields a strict inequality at the beginning of \eqref{eq:chain} in Step 2, which gives a contradiction. Hence,  conclusion (1) does not hold and only conclusion (2) can occur
in these dimensions.
 \end{proof}

\section*{AI disclosure}

The authors acknowledge the use of ChatGPT to check the correctness of some proofs, to improve the presentation of the manuscript, and to assist in finding the example presented in Section \ref{sec:domains condition C}. All proof strategies and mathematical arguments were developed by the authors. The authors take full responsibility for the entire content of the paper.


\begin{thebibliography}{99}

\bibitem[Aub76]{Aubin_1976}
T.~Aubin, \textit{Problèmes isopérimétriques et espaces de Sobolev}, J. Differ. Geom. \textbf{11} (1976), no.~4, 573--598.

\bibitem[Bre11]{Brezis_2011}
H.~Brezis, \textit{Functional Analysis, Sobolev Spaces and Partial Differential Equations}, Springer, New York, NY, 2011.

\bibitem[BL83]{BL83}
H.~Brezis and E.~Lieb,
\textit{A Relation Between Pointwise Convergence of Functions and Convergence of Functionals},
Proc. Am. Math. Soc.
\textbf{88} (1983), no. 3, 486--490.

\bibitem[BN83]{BN}
H.~Brezis, L.~Nirenberg, \textit{Positive solutions of nonlinear elliptic equations involving critical sobolev exponents}, Comm. Pure Appl. Math. \textbf{36} (1983), no.~4, 437--477.

\bibitem[CGS89]{CaGiSp}
L.~A.~Caffarelli, B.~Gidas, J.~Spruck, 
\textit{Asymptotic symmetry and local behavior of semilinear elliptic equations with critical Sobolev growth}, 
Commun. Pure Appl. Math. \textbf{42} (1989), no. 3, 271--297.

\bibitem[Cha81]{Chang_1981}
K.-C.~Chang, \textit{Variational methods for non-differentiable functionals and their applications to partial differential equations}, J. Math. Anal. Appl. \textbf{80} (1981), no.~1, 102--129.

\bibitem[CZ12]{Chen_Zou_2012}
Z.~Chen, W.~Zou, \textit{Positive Least Energy Solutions and Phase Separation for Coupled Schrödinger Equations with Critical Exponent}, Arch. Ration. Mech. Anal. \textbf{205} (2012), no.~2, 515--551.

\bibitem[CZ15]{Chen_Zou_2015}
Z.~Chen, W.~Zou, \textit{Positive least energy solutions and phase separation for coupled Schrödinger equations with critical exponent: higher dimensional case}, Calc. Var. Partial Differ. Equ. \textbf{52} (2015), no.~1--2, 423--467.

\bibitem[CSS26]{Clapp_Saldaña_Schiera_2026}
M.~Clapp, A.~Saldaña, D.~Schiera, \textit{Sublinear elliptic equations with a sharp change of sign in the nonlinearity}, preprint (2026), arXiv:2603.11766.

\bibitem[Cla76]{Clarke_1976}
F.~H.~Clarke, \textit{A New Approach to Lagrange Multipliers}, Math. Oper. Res. \textbf{1} (1976), no.~2, 165--174.

\bibitem[CCO09]{CCO}
M.~Colin, T.~Colin, M.~Ohta, \textit{Stability of solitary waves for a system of nonlinear Schrödinger equations with three wave interaction}, Ann. Inst. Henri
Poincar\'e, Anal. Non Lin\'eaire \textbf{26} (2009), no.~6, 2211--2226.

\bibitem[FGS26]{FaGiSo26}
H.~Fan, L.~Giaretto, N.~Soave, \textit{On least energy solutions for a critical nonlinear Schr\"odinger system with $K$-wise interaction: the cooperative case}, Preprint (2026).

\bibitem[Gia26]{Gi}
	L.~Giaretto, \emph{On elliptic systems with $k$-wise interactions in the strong competition regime: uniform H\"older bounds and properties of the limiting configurations}, preprint, arXiv:2603.10949 (2026).

\bibitem[GS26]{Giaretto_Soave_2026}
L.~Giaretto, N.~Soave, \textit{On least energy solutions for a nonlinear Schrödinger system with K-wise interaction}, Nonlinear Anal. \textbf{262} (2026), 113938.

\bibitem[LYZ23]{LYZ23}
T.~Liu, S.~You, W.~Zou, 
\textit{Least energy positive solutions for $d$-coupled Schr\"odinger systems with critical exponent in dimension three}. 
J. Differ. Equ. \textbf{367} (2023), 40--78.

\bibitem[LQZ24]{LQZ} 
C.~Luo, Y.~Qiang, D.~Zwicker, \textit{Beyond pairwise: Higher-order physical interactions affect phase separation in multicomponent liquids}, Phys. Rev. Res. \textbf{6}
(2024), 033002.

\bibitem[NT+10]{NoTaTeVe}
	B.~Noris, H.~Tavares, S.~Terracini, G.~Verzini, \textit{Uniform H\"older bounds for nonlinear Schr\"odinger systems with strong competition}, Commun. Pure Appl. Math. \textbf{63} (2010), no.~3, 267--302.

\bibitem[Pet14]{Pet}
D.~S.~Petrov, \textit{Three-body interacting bosons in free space}, Phys. Rev. Lett. \textbf{112} (2014), 103201.

\bibitem[Pom10]{Pom}
A.~Pomponio, \textit{Ground states for a system of nonlinear Schr\"odinger equations with three wave interaction}, J. Math. Phys. \textbf{51} (2010), no.~9, 093513.

\bibitem[ST+16]{STTZ16}
N. Soave, H. Tavares, S. Terracini, and A. Zilio.
\textit{H\"{o}lder bounds and regularity of emerging free boundaries for
strongly competing {S}chr\"{o}dinger equations with nontrivial grouping},
 Nonlinear Anal. \textbf{138} (2016), 388--427.

\bibitem[ST24a]{SoTe1}
	N.~Soave, S.~Terracini, \textit{On some singularly perturbed elliptic systems modeling partial segregation: uniform H\"older estimates and basic properties of the limits}, preprint, arXiv:2409.11976 (2024).

\bibitem[ST24b]{SoTe2}
	N.~Soave, S.~Terracini, \textit{On partially segregated harmonic maps: optimal regularity and structure of the free boundary}, preprint, arXiv:2410.23976 (2024).

\bibitem[SZ15]{SZ15}
N.~Soave and A.~Zilio.
\textit{Uniform bounds for strongly competing systems: the optimal {L}ipschitz case},
 Arch. Ration. Mech. Anal. \textbf{218} (2015), no.~2, 647--697. 

\bibitem[Tal76]{Talenti_1976}
G.~Talenti, \textit{Best constant in Sobolev inequality}, Ann. Mat. Pura Appl. \textbf{110} (1976), no.~1, 353--372.

\bibitem[TY20]{TY20}
H.~Tavares, S.~You, \textit{Existence of least energy positive solutions to Schr\"odinger systems with mixed competition and cooperation terms: the critical case}. 
Calc. Var. Partial Differ. Equ. \textbf{59} (2020), no.~1, Paper No. 26, 35 p.

\bibitem[TYZ22]{TYZ22}
H.~Tavares, S.~You, W.~Zou, 
\textit{Least energy positive solutions of critical Schr\"odinger systems with mixed competition and cooperation terms: the higher dimensional case}, 
J. Funct. Anal. \textbf{283} (2022), no.~2, Article ID 109497, 50 p.

\bibitem[ZZ26]{ZZ26}
Z. Zhang, Z. Zhang, 
\textit{Uniform global Lipschitz bounds of solutions to nonlinear Schr\"odinger systems with strong competition},
J. Funct. Anal. \textbf{291} (2026), no.~3, 42 p.

\end{thebibliography}
\end{document}